\documentclass[3p]{elsarticle}

\usepackage{lineno,hyperref}
\usepackage{amsmath}
\usepackage[dvipsnames]{xcolor}
\usepackage{listings}
\usepackage{multirow}
\usepackage{floatrow}
\usepackage{amssymb}

\newtheorem{definition}{Definition}[section]

\usepackage{mathtools}
\DeclarePairedDelimiterX{\norm}[1]{\lVert}{\rVert}{#1}

\usepackage{cleveref}
\usepackage{hyperref}

\usepackage[symbol]{footmisc}

\biboptions{comma,square,sort&compress}

\begin{document}

\begin{frontmatter}

\title{{\color{black}High-performance \texttt{dune} modules for solving large-scale, strongly anisotropic elliptic problems with applications to aerospace composites}}

\author[be]{R. Butler}
\address[be]{Department of Mechanical Engineering, University of Bath, UK.}
\author[e,Tur]{T. Dodwell }
\address[e]{College of Engineering, Mathematics and Physical Sciences, University of Exeter, UK.}
\address[Tur]{The Alan Turing Institute, London, NW1 2DB, UK.}
\author[m]{A. Reinarz\footnotemark}
\address[m]{Institute of Informatics, Technical University of Munich, Germany.}
\author[e]{A. Sandhu}
\author[h,bm]{R. Scheichl}
\address[h]{Institute for Scientific Computing, University of Heidelberg, Germany.}
\address[bm]{Department of Mathematical Sciences, University of Bath, UK.}
\author[h]{L. Seelinger}

\footnotetext{Corresponding Author: Anne Reinarz \texttt{reinarz@in.tum.de}}

\begin{abstract}
{\color{black}
The key innovation in this paper is an open-source, high-performance iterative solver for high contrast, strongly anisotropic elliptic partial differential equations implemented within \texttt{dune-pdelab}. The iterative solver exploits a robust, scalable two-level additive Schwarz preconditioner, GenEO (Spillane et al. 2014). The development of this solver has been motivated by the need to overcome the limitations of commercially available  modeling tools for solving structural analysis simulations in aerospace composite applications. Our software toolbox \texttt{dune-composites} encapsulates the mathematical complexities of the underlying packages within an efficient C++ framework, providing an application interface to our new high-performance solver. We illustrate its use on a range of industrially motivated examples, which should enable other scientists to build on and extend \texttt{dune-composites} and the GenEO preconditioner for use in their own applications. We demonstrate the scalability of the solver on more than 15,000 cores of the UK national supercomputer \textsc{Archer}, solving an aerospace composite problem with over 200 million degrees of freedom in a few minutes. This scale of computation brings composites problems that would otherwise be unthinkable into the feasible range. To demonstrate the wider applicability of the new solver, we also confirm the robustness and scalability of the solver on SPE10, a challenging  benchmark in subsurface flow/reservoir simulation.}  
\end{abstract}


\end{frontmatter}



\noindent{\bf PROGRAM SUMMARY}\\

\begin{small}
\noindent
{\em Program Title:} \texttt{dune-composites}                                      \\
{\em Licensing provisions:} BSD 3-clause                                    \\
{\em Programming language:} C++                                  \\


  \noindent{\em Nature of problem:}\\
\texttt{dune-composites} is designed to solve anisotropic linear elasticity equations for anisotropic, heterogeneous materials, e.g. composite materials. To achieve this, our contribution also implements a new preconditioner in \texttt{dune-pdelab}. 

  \noindent{\em Solution method:}\\
The anisotropic elliptic partial differential equations are solved via the finite element method. The resulting linear system is solved via an iterative solver with a robust, scalable two-level overlapping Schwarz preconditioner: GenEO.
 

\end{small}


\section{Introduction}

{\color{black}
\noindent Across the physical sciences, elliptic partial differential equations (PDEs) naturally arise as mathematical models of the equilibrium state of a system. Classical examples include the distribution of temperature in a body, the flow of fluid in a porous medium and, the particular application of interest in this paper, the equilibrium of forces acting on a material. The most widely used approach to solve such PDEs is the finite element method~\cite{Hug00}, which results in a sparse system of equations. For systems which exhibit multiple scales and large spatial variations in model parameters the system of equations can be very ill-conditioned and extremely large, e.g.~contain $>10^7$ degrees of freedom. The design of robust and scalable solvers that do not require laborious tuning and are capable of exploiting the power of modern distributed computers, is essential. In this paper, we describe the design and implementation of a robust two-level additive Schwarz preconditioner for parallel Krylov based iterative solvers in \texttt{dune-pdelab} \cite{dune}. By also creating a bespoke module for analysis of composite structures \texttt{dune-composites}, we demonstrate the capabilities of this new solver on industrially motivated aerospace composite problems with over 200 million degrees of freedom. To also demonstrate the wider applicability of the new solver, we demonstrate the robustness and scalability of the solver on the challenging, classical SPE10 benchmark \cite{spe10,eike} in subsurface flow/reservoir simulation.}

\subsection{Motivating Computational Challenge in Aerospace Composites}

\smallskip
\noindent 
Scientific advances in aerospace composite design and materials offer exciting engineering opportunities, making them the material of choice for many modern aircraft (e.g Airbus A350, Boeing 787). However, composite manufacturers face huge challenges in designing and making complex components quickly enough to remain commercially competitive.
There is a growing realisation in both academia and industry that to meet ambitious global growth targets, composites manufacturers should {\em `reduce time, cost and risk to market through the use of validated simulation packages'} \cite{CompositeStrategy}. Currently simulation capabilities allowing high-fidelity full-scale analysis of a composite structure are not openly available. But, why is this? What makes this analysis of large scale composite structure so challenging?

\begin{figure}[ht!]
\includegraphics[width=0.9\linewidth]{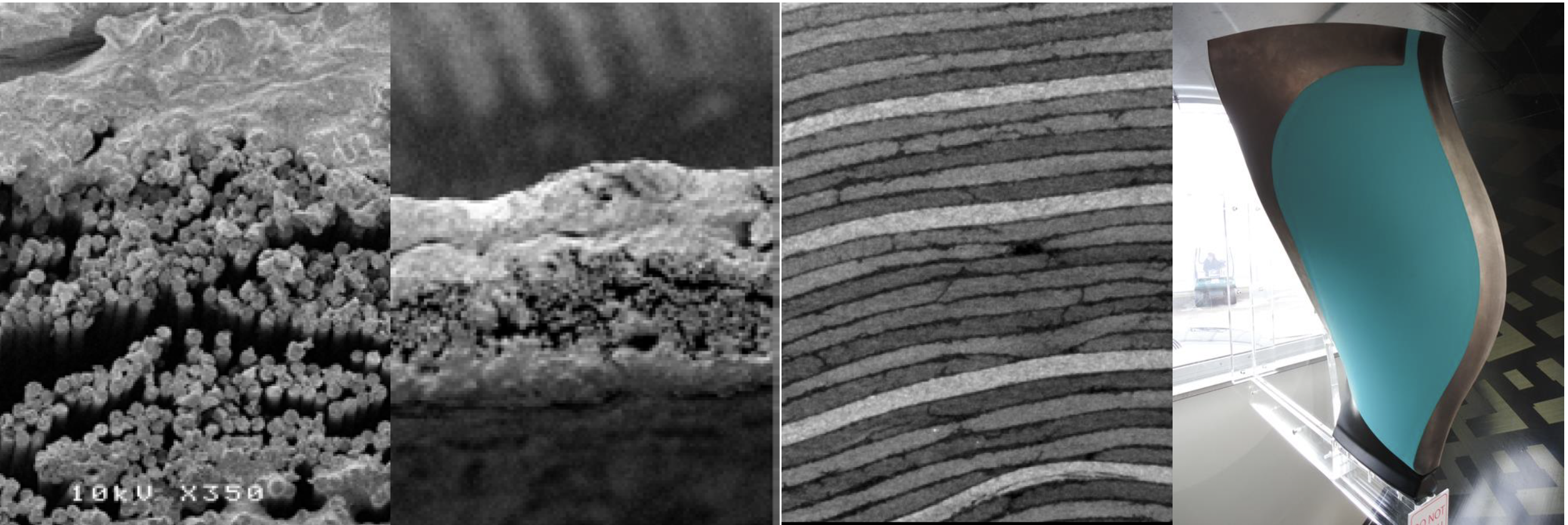}
\caption{\label{fig:intro}The range of scales introduces significant complexity in the analysis of aerospace composites. From Left to Right:
fibre/resin scale ($5\mu$m), ply scale (0.25mm), laminate scale (5-30mm) to the structure ($>1$m)(Composite Fan Blade, Right).}
\end{figure}

\smallskip\noindent
When we apply classical finite element (FE) analysis to a composite structure the problem reduces to finding a vector of displacements ${\bf u}^{(i)}\in \mathbb{R}^3$ at each of the $N$ nodes within a FE mesh. This leads to the sparse system of FE equations \cite{Hug00}:
\begin{equation}\label{eqn:MatrixEq}
{\bf A} \tilde{{\bf u}} = {\bf b}, \quad \mbox{where} \quad  \tilde{{\bf u}} = [{\bf u}_h^{(1)}, \ldots, {\bf u}_h^{(N)}]^T,
\end{equation}
$\bf{A}$ is the global stiffness matrix and $\bf{b}$ is the load vector arising from the applied boundary conditions or loading. In solving the linear system \eqref{eqn:MatrixEq}, we face two significant mathematical challenges:
\begin{itemize}
\item {\bf Scale of Calculations.} 
Composite materials are manufactured from thin fibrous layers, less than $1$mm thick, separated by even thinner resin interfaces of thickness less than $0.05$mm, yet entire component parts are generally several meters long, Fig. \ref{fig:intro}. To resolve stresses and accurately predict failure, several elements need to be placed through each layer \cite{dunecomp}. Naturally this means that the number of nodes $N$ is very large. As an example in this paper, we model a $1$m section of a wing box given in Section \ref{sec:wingbox}, while resolving the resin interfaces, giving in total $200$ million degrees of freedom. Solving linear systems of this size requires specialised, parallel solvers. Current industry standard tools, such as \textsc{Abaqus} \cite{Abaqus}, are not able to deal with these problem sizes, largely due to limitations of the parallel solvers employed.

\item {\bf Material Anisotropy.} Central to the benefits of composite structures is the inclusion of directional fibres, which gives them an excellent weight to stiffness ratio under a particular loading. This means, there is a large contrast ($\sim 1:40$) in mechanical properties within a single layer of composite, related to the fibre direction(s) and those directions dominated by the stiffness of the matrix material (typically a toughened epoxy resin). In the FE discretisation, this leads to a stronger coupling between degrees of freedom in the fibre direction, as opposed to those in the orthogonal directions. 
The fibrous layers are stacked with different fibre orientations to form a laminate, adding an additional level of complexity. The fibre directions act as stiff constraints on the deformation, whilst the weak connections give rise to low-energy mechanics within the structure. Mathematically, this causes significant numerical challenges in solving \cref{eqn:MatrixEq} via iterative solvers, since the system is very ill-conditioned. For such cases, classical iterative solvers (required to address Challenge 1) converge very slowly.
\end{itemize}

\smallskip\noindent
The usual approach to tackle both these challenges is to apply the parallel iterative solvers to a {\em  preconditioned} version of \cref{eqn:MatrixEq}. The task then is to develop an operator $\textbf{M}^{-1}$, which is computationally cheap to construct, such that $\textbf{M}^{-1}\textbf{A}\tilde{{\bf u}} = \textbf{M}^{-1}{\bf b}$ is better conditioned. The most widely used preconditioners for iterative solvers for \eqref{eqn:MatrixEq} in both commercial and scientific FE codes are Algebraic Multigrid (AMG) methods \cite{Bla07,boomerAMG}. They have demonstrated excellent scalability for a broad class of problems over thousands of processors, and have the advantage of working only on the matrix equations \eqref{eqn:MatrixEq}, so that they can be applied `black-box'. As a preconditioner, AMG constructs the matrix ${\bf M}$ by repeatedly coarsening the full matrix ${\bf A}$ through recursive aggregation over the degrees of freedom . The aggregation process is algebraic and based on the fact that the solution at two neighbouring nodes will be similar if they are `strongly connected'. The success of an AMG preconditioner depends on this aggregation process. As discussed above, the connectivity of degrees of freedom within a laminate is very complex even for a simple laminate and in our numerical experiments the performance of all AMG preconditioners that we tested was prohibitively poor. In particular, this includes off-the-shelf AMG used in the commercial software \textsc{Abaqus} \cite{Abaqus}, as well as tailored aggregation strategies in the AMG preconditioners provided through \texttt{dune-istl} \cite{Bla07} and \texttt{Hypre} \cite{boomerAMG}.

\smallskip\noindent
This initial testing of AMG preconditioners highlighted the need for the development and implementation of a robust preconditioner with respect to both problem size, material heterogeneity and anisotropy for large-scale composite based applications. Over the last decade there has been significant effort from the domain decomposition community to develop scalable and robust preconditioners suitable for parallel computation. One such preconditioning approach is provided by the additive Schwarz framework \cite{Tos04}. The domain is decomposed into overlapping subdomains, which in our case each correspond to one processor, and the subdomain's local stiffness matrix is inverted on each processor using a direct solver. This ``one level'' approach is not sufficient for very large problems and global information in the form of a coarse space must be added. In \texttt{dune-composites}, we use GenEO  \cite{Spi14} to construct a coarse space by combining low energy eigenvectors of the local subdomain stiffness matrices using a partition of unity. The resulting preconditioner leads to an almost optimal scaling with respect to problem size and number of processors, allowing us to successfully tackle large industrially important problems with over 200 million degrees of freedom.

\subsection{The Contributions of this Paper}

\smallskip\noindent
\texttt{dune-composites} is a high-performance composite FE package built on top of \textsc{Dune} (Distributed and Unified Numerics Environment), an open source modular toolbox for solving partial differential equations (PDEs) with grid-based methods \cite{Bas08a,Bas08b,Bas08c}. Based on the core \textsc{Dune} philosophy, \texttt{dune-composites} is written using C++ and exploits modern inheritance and templating programming paradigms. {\color{black} It is open-source and publicly available at \url{https://dune-project.org/modules/dune-composites/}.} The package provides a codebase with the following key features:

\begin{itemize}
\item implementation and interface to a novel, robust preconditioner called \textsc{GenEO} \cite{dunecomp,Spi14} for parallel Krylov solvers, which exhibits excellent scalability over thousands of processors on Archer, the UK national HPC system. {\color{black} Since release~2.6, the preconditioner is provided as part of \texttt{dune-pdelab}, available at \url{https://dune-project.org/modules/dune-pdelab/} (initially it had been developed within the \texttt{dune-composites} module)};

\item interfaces to handle composite applications, including stacking sequences, complex part geometries, defects and non-standard boundary conditions, such as multi-point constraints or periodicity;

\item to overcome shear locking of standard FEs, mesh stabilisation strategies to support reduced integration \cite{Bel84}, as well as a new 20-node 3D serendipity element (with full integration) have been implemented;

\item interfaces to other state-of-the-art parallel solvers (\& preconditioners) in \texttt{dune-istl} \cite{Bla07} and \texttt{Hypre} \cite{boomerAMG};

\item a code structure which supports both engineering end-users, and those requiring flexibility to extend any aspect of the code in a modular way to introduce new applications, solvers or material models.

\end{itemize}

\smallskip\noindent
The purpose of this paper is to highlight the novel mathematical aspects of the code and document its structure. We illustrate its use through a range of industry motivated examples, which should enable other scientists to build on and extend \textit{dune-composites} for use in their own applications. We begin by outlining the mathematical formulation of the new robust preconditioner and its implementation on a distributed memory computer in Sec.~\ref{sec:solver}. We then provide details of the structure and salient features of the code in Sec.~\ref{sec:code}. Finally in Sec.~\ref{sec:examples}, through the use of a series of example problems, we provide details of how to implement, build and run your own applications. We also use these examples as an opportunity to demonstrate the computational efficiency of \texttt{dune-composites}.

\section{Preliminaries : Anisotropic Elasticity Equations and their Finite Element Discretisation}

\smallskip\noindent
A composite structure occupies the domain $\Omega \subset \mathbb R^3$ with the boundary $\Gamma$ and a unit, outward normal vector ${\bf n} \in \mathbb R^3$. {\color{black} At each point ${\bf x}\in \Omega$ we define a vector-valued displacement ${\bf u}({\bf x}): \Omega \rightarrow \mathbb R^3$.  In each of these three global directions }the boundary may contain a Dirichlet component $\Gamma_D^{(i)}$ and a Neumann component $\Gamma_N^{(i)}$, such that
\begin{equation}\label{eq-boundary}
\Gamma = \Gamma_D^{(i)} \cup \Gamma_N^{(i)} \quad \mbox{and} \quad \Gamma_D^{(i)} {\color{black}\cap }\Gamma_N^{(i)} = \emptyset, \quad i = x, y, z.
\end{equation}

\smallskip\noindent Let $\sigma_{ij}$ denote the Cauchy stress tensor and ${\bf f}({\bf x}): \Omega \rightarrow \mathbb R^3$ the body force per unit volume. The infinitesimal strain tensor, is defined as the symmetric part of the displacement gradients
\begin{equation}\label{eq-strain}
\epsilon_{ij}({\bf u}) = \frac{1}{2}\left(u_{i,j} + u_{j,i}\right),
\end{equation}
{\color{black} where $u_{i,j} = \frac{\partial u_i}{\partial x_j} $.}
The strain tensor is connected to Cauchy stress tensor via the generalised Hooke's law
\begin{equation}\label{eq-hooke}
\sigma_{ij}({\bf u}) = C_{ijkl}({\bf x})\epsilon_{kl}({\bf u}).
\end{equation}
$C_{ijkl}({\bf x})$ is a symmetric, positive definite fourth order tensor. A composite laminate is made up of a stack of composite layers (or plies $\sim 0.2$mm), separated by a very thin layer of resin ($15\mu$m). A single composite layer is modelled as a homogeneous orthotropic elastic material, characterised in general by $9$ parameters and {\color{black} a vector of  orientations ${\bf\theta}$}. Resin interfaces are assumed isotropic, defined by just $2$ scalar {\color{black}(Lam\'e)} parameters. {\color{black} These fibres are aligned in local coordinates and can be rotated in any direction using standard tensor rotations, for more details see e.g. \cite{Ting92}.}

\smallskip\noindent
Given functions $h_i: \Gamma_D^{(i)} \rightarrow \mathbb R$ and $g_i: \Gamma_N^{(i)} \rightarrow \mathbb R$, prescribing the Dirichlet and Neumann boundary data (for each component), we seek the unknown displacement field ${\bf u}({\bf x})$, which satisfies the force equilibrium equations and the boundary conditions,
\begin{equation}\label{eq:strongForm}
\nabla \cdot {\underline{\sigma} (\bf u)} + {\bf f} = 0, \quad  {\bf x} \in \Omega, \quad
u_i = h_i\; \mbox{for} \quad {\bf x} \in \Gamma_D^{(i)} \quad \mbox{and} \quad \sigma_{ij}n_j = g_i \quad \mbox{for} \; {\bf x} \in \Gamma_N^{(i)},
\end{equation}
as well as eqs.~\eqref{eq-strain} and \eqref{eq-hooke}. Then, we define the function space for each component of displacement $u_i$ to be
\begin{equation}\label{eq-fctspace}
V^{(i)} := \{v \in H^1(\Omega) : v_i({\bf x}) = h_i \;, \ {\bf x} \in \Gamma_D^{(i)} \} ,
\end{equation}
\noindent leading to the weak formulation of \eqref{eq:strongForm}, of finding ${\bf u} \in V := V^{(1)} \otimes V^{(2)} \otimes V^{(3)}$ such that
\begin{equation}\label{eq-fe}
a({\bf u}, {\bf v}) := \int_\Omega \sigma_{ij}({\bf u})\epsilon_{ij}({\bf v}) \; dx = \int_{\partial\Omega_n} \sigma_{ij}n_jv_i\;ds - \int_\Omega f_iv_i\;dx := b({\bf v}), \quad \forall \mathbf{v} \in V.
\end{equation}

\noindent
We consider the discretisation of the variational equations \eqref{eq-fe} with  conforming FEs on a mesh $\mathcal T_h$ on $\Omega$. Let $V_h \subset V$ denote the restriction of $V$ onto a FE space on $\mathcal T_h$ and seek an approximation ${\bf u}_h \in V_h$ such that
\begin{equation}\label{eqn:bilinear}
a({\bf u}_h,{\bf v}_h) - b({\bf v}_h) = 0, \quad \mbox{for all} \quad {\bf v}_h \in V_h.
\end{equation}
We block together displacements from all three space dimensions, so that ${\bf u}^{(i)}_h \in \mathbb B = \mathbb R^3$ denotes the vector of displacement coefficients containing all space components associated with the $i^{th}$ basis function. We introduce the (vector-valued) FE basis for $V_h$ defined by the spanning set of (vector-valued) shape functions $\{\boldsymbol\phi^{(i)}({\bf x})\}_{i=1}^N$. These are the normal scalar shape functions, repeated for each displacement component. Therefore the vector displacement at a point is given by ${\bf u}_h ({\bf x}) = \sum_{i=1}^N ({\bf u}_h^{(i)})^T \;\boldsymbol \phi^{(i)}({\bf x})$. The choice of basis converts \eqref{eqn:bilinear} into a symmetric positive-definite (spd) system of algebraic equations
\begin{equation}\label{eqn:fe_matrix_system}
{\bf A}{\bf \tilde{u}} = {\bf b} \quad \mbox{where} \quad {\bf A} \in \mathbb B^N \times \mathbb B^N \quad \mbox{and} \quad {\bf b} \in \mathbb B^N
\end{equation}
where the blocks in the  global stiffness matrix and in the load vector, for any $i, j = 1, \ldots, N$, are given by ${\bf A}_{ij} = a(\boldsymbol \phi^{(i)}, \boldsymbol \phi^{(j)})$ and ${\bf b}_i = b(\boldsymbol \phi_i)$. The vector $\tilde{{\bf u}} = [{\bf u}_h^{(1)}, \ldots, {\bf u}_h^{(N)}]^T \in \mathbb B^N$ is the block vector of unknown FE coefficients.
System \eqref{eqn:fe_matrix_system} is assembled element-wise from \eqref{eq-fe}, using Gaussian integration.

\section{A Robust, Scalable, Parallel Iterative Solver for Composite Structures}\label{sec:solver}

\noindent 
The key innovation of {\em dune-composites}, as a software package, is the design and implementation of a highly robust, scalable parallel iterative solver for composite applications. {\color{black} This solver is applicable to a more general class of problems and has been made available in the \textit{dune-pdelab} module.} In this section, we provide the mathematical and implementation details of the new solver. Apart from new types of FEs that had to be implemented in {\sc Dune}, the remainder of the package largely provides interfaces to handle the set-up for complex composite applications.

\subsection{Krylov Subspace Methods Preconditioned with Two-Level Additive Schwarz Methods}

\noindent In \texttt{dune-composites} we use {\em preconditioned Krylov subspace methods} both in sequential and parallel, as provided by {\sc Dune}'s ``Iterative Solver Template Library'' \texttt{dune-istl} \cite{Bla07}. Krylov subspace methods are iterative solvers which construct a sequence of approximations ${\bf u}^{(k)}$ in the $k$-dimensional subspace:
\begin{equation}\label{eqn:KrylovSubspace}
\mathcal K_k = \mbox{span}\{ {\bf r}, {\bf A r}, {\bf A}^2 {\bf r}, \ldots,{\bf A}^{k-1} {\bf r} \} \subset \mathbb R^n
\end{equation}
where ${\bf r} =  {\bf b} - {\bf A} {\bf \tilde u}^{(0)}$ is the initial residual. The simplest Krylov subspace method for a symmetric positive-definite matrix ${\bf A}$ is the Conjugate Gradient method (CG), first introduced by Hestenes and Stiefel (1952). In each step, the approximate solution ${\bf \tilde{u}}^{(k)}$ is updated by adding the search direction ${\bf d}^{(k)}$ scaled by a factor chosen to minimise the energy norm over the space ${\mathcal K_k}$. The search directions are chosen to be ${\bf A}$-orthogonal to all previous direction i.e. $\langle {\bf d}^{(k)}, {\bf A}{\bf d}^{(k')}\rangle = 0$ for $k^\prime < k$. The method iterates until the residual norm (or ``energy'') $\|{\bf r}^{(k)}\|$ reduces below a user defined tolerance. Importantly, the convergence rate of CG depends on the spectral properties of the matrix ${\bf A}$, see e.g.~\cite{Saa03}. In particular, it can be bounded proportionally to the square root of the condition number $\kappa$, defined as the ratio between its largest and smallest eigenvalue. A large value, as usually seen in composites, indicates that the system ${\bf A}{\tilde {\bf u} } = {\bf b}$ is ill-conditioned. This means that ${\bf u}$ is very sensitive to small changes in ${\bf b}$. For such cases, iterative solvers converge very slowly or even not at all, particularly when the problem size increases. A remedy is to {\em precondition} the system, that is to develop an operation ${\bf M}^{-1}$ which is computationally cheap to construct and apply (in parallel) such that ${\bf M}^{-1}{\bf A}\tilde{{\bf u}} = {\bf M}^{-1}{\bf b}$ is better conditioned and CG solvers converge quickly.

\begin{figure}
\includegraphics[width=0.28\linewidth]{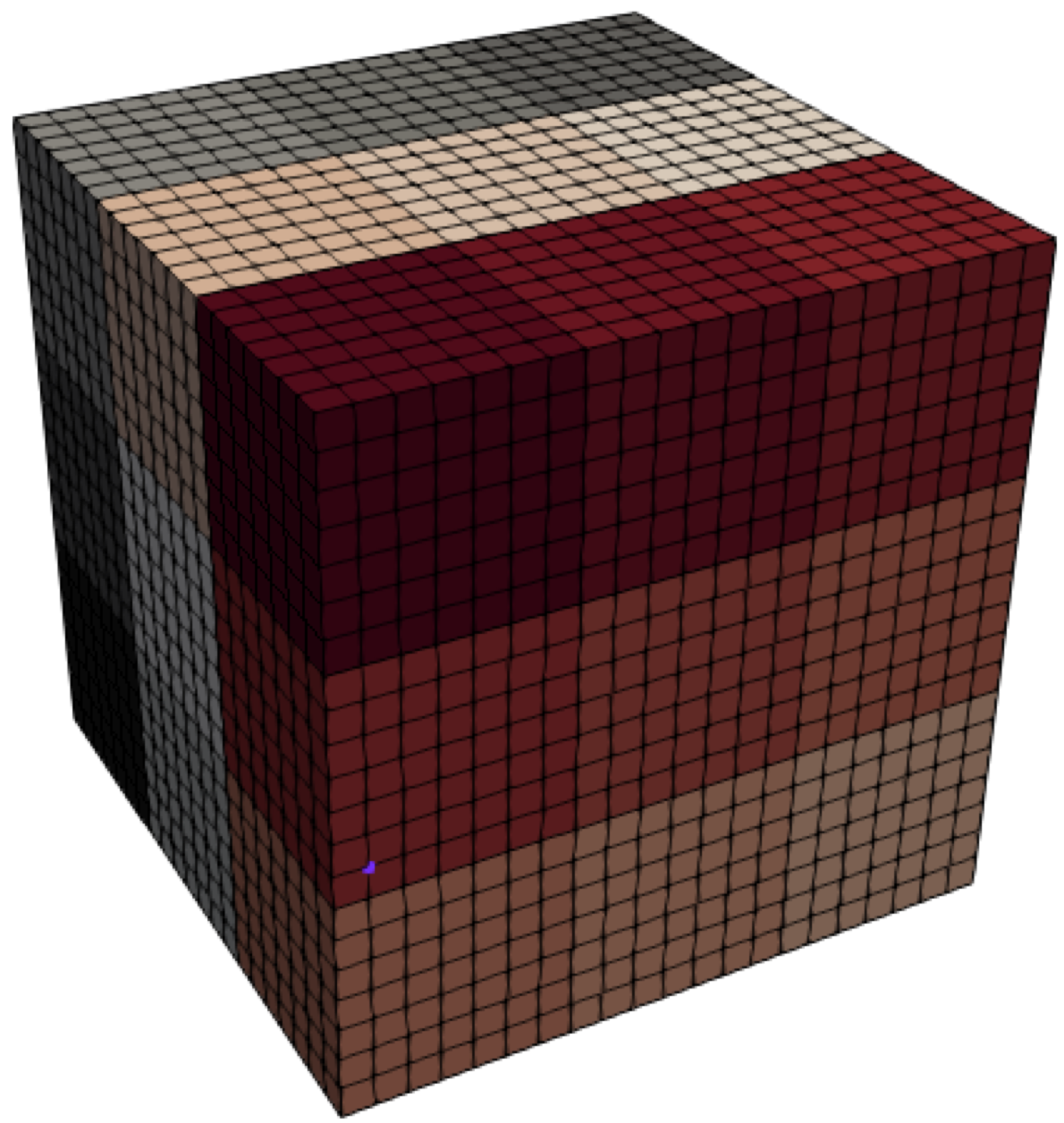}\includegraphics[width=0.28\linewidth]{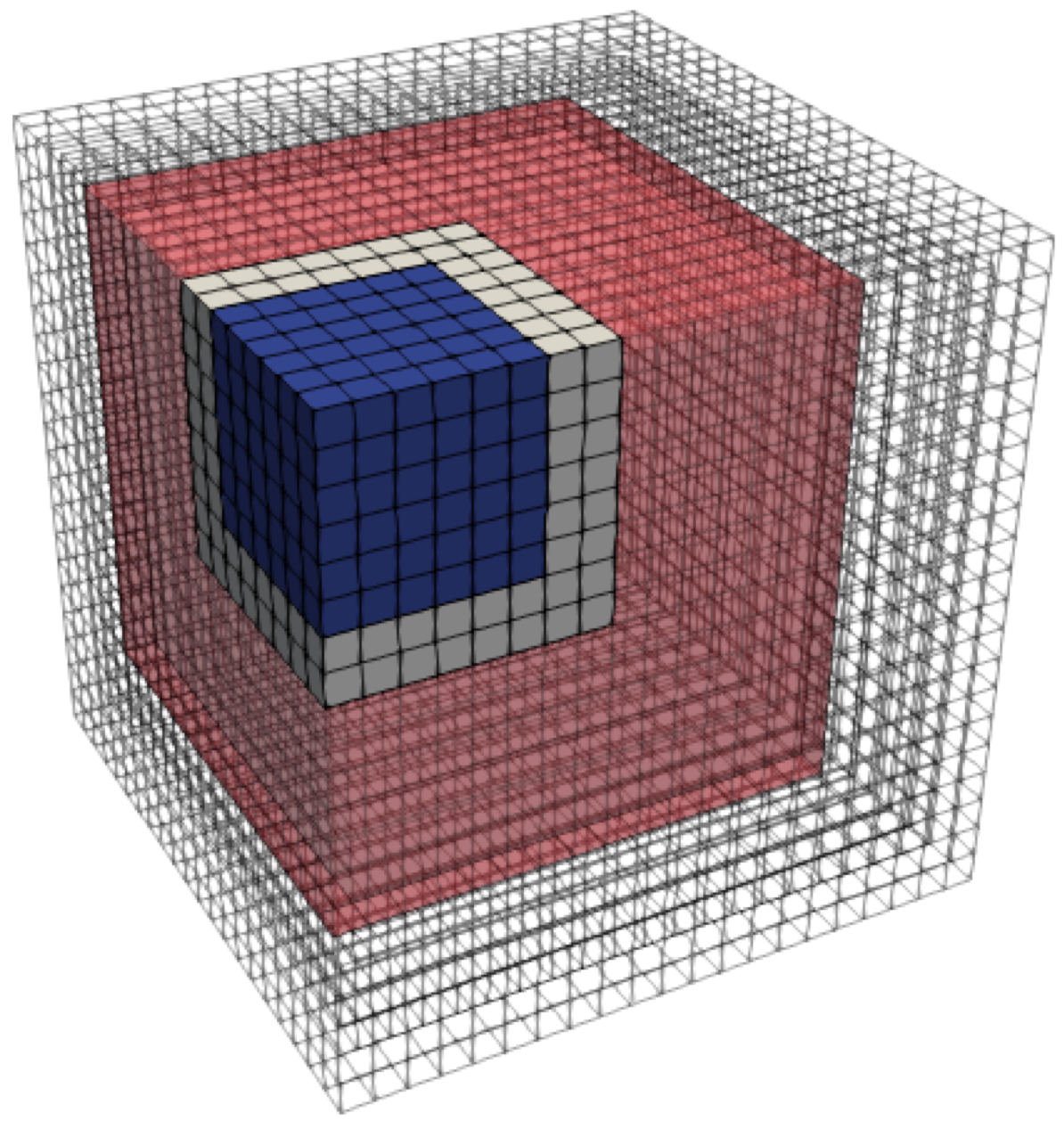}\includegraphics[width=0.35\linewidth]{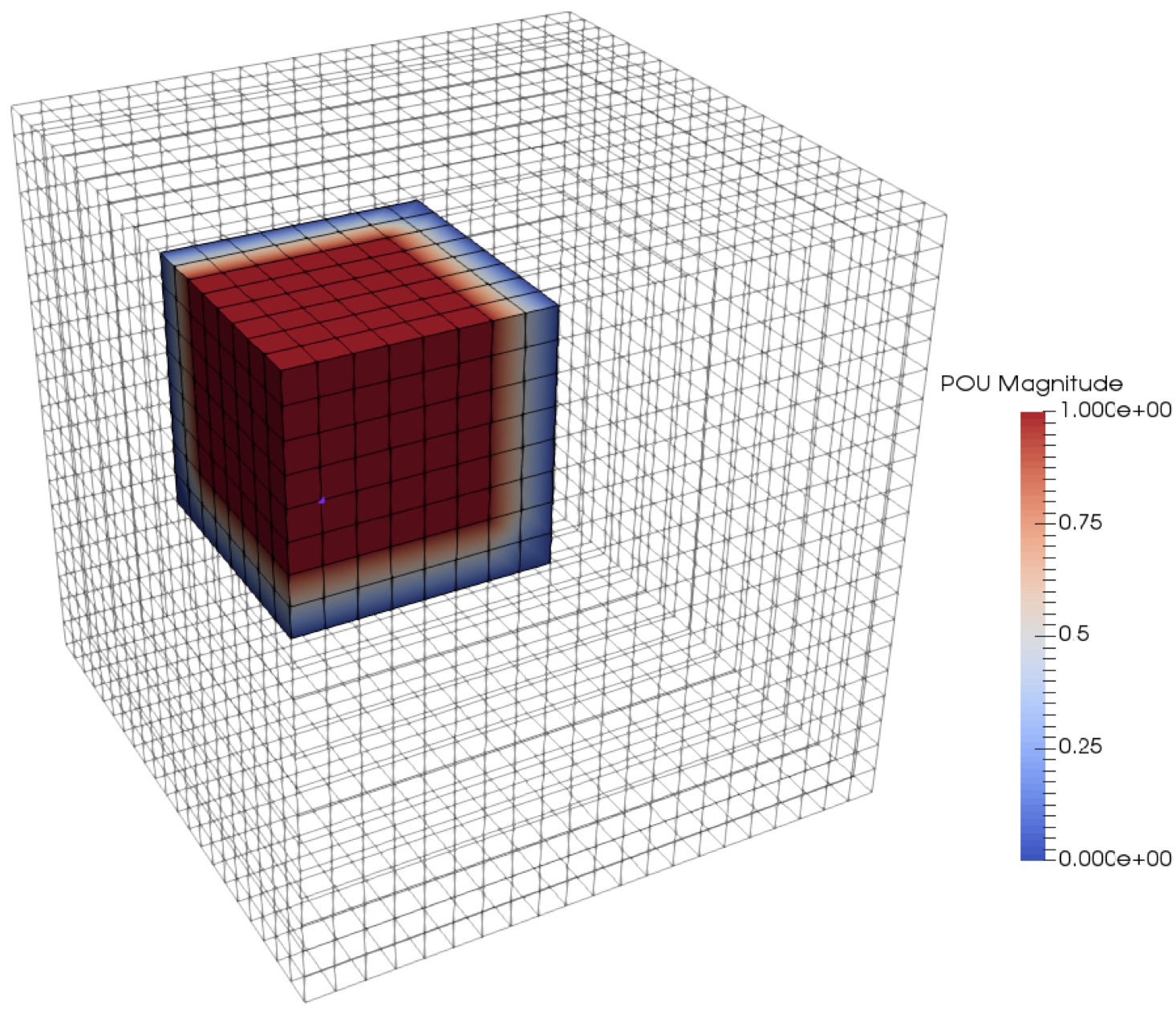}
\caption{\label{fig:domaindecomp}(Left) Domain $\Omega$ partitioned into non-overlapping subdomains $\Omega_j'$ where colouring differentiates independent subdomains. (Middle) Shows overlapping subdomain $\Omega_j$ with a single layer of overlap ($O=1$). Overlap region $\Omega^\circ_j$ is shown in white. Transparent red regions show cells of the grid which belong to 'nearest neighbour' processors. (Right) Shows partition of unity (PoU) operator $\Xi_j$ on a single processor, defined as in \eqref{eqn:pou}.}
\end{figure}

\smallskip\noindent
In \texttt{dune-composites} our main preconditioner is a two level additive Schwarz method. 
To construct this method we partition our domain $\Omega$ into a set of non-overlapping subdomains $\Omega_j^{'}$ for $j=1$ to $N$ resolved by $\mathcal T_h$, as shown in Fig. \ref{fig:domaindecomp} (left). Each subdomain $\Omega_j^{'}$ is extended by $O$-layers of elements to give the overlapping subdomains $\Omega_j$, Fig. \ref{fig:domaindecomp} (middle).
For each subdomain $1 \leqslant j \leqslant N$, we denote the restriction of $V_h$ to $\Omega_j$ by $V_h ( \Omega_j)$, whilst the space of FE functions with support contained in $\Omega_j$ is called $V_{h, 0} ( \Omega_j)$. 

\smallskip\noindent
{\em Remark:} In \texttt{dune-composites} the user can define the initial non-overlapping decomposition (or a default is used), the overlapping process is handled by {\sc Dune}'s parallel structured grid class \texttt{Dune::YaspGrid} \cite{Bas08a}.

\smallskip\noindent
Any function $v \in V_{h, 0} ( \Omega_j)$  is mapped onto $V_h$ by the prolongation operator $R_j^T : V_{h, 0} ( \Omega_j) \rightarrow V_h$, which extends $v$ by zeros, so that
 \begin{displaymath}
    R_j^T v({\bf x}) = \begin{cases}v({\bf x}), & {\bf x} \in  \Omega_j \\  0,  & {\bf x} \in \Omega\backslash \Omega_j\end{cases}.
 \end{displaymath}
We therefore note that the restriction operator $R_j:V_h \rightarrow V_{h, 0} ( \Omega_j)$. In matrix form the restriction and prolongation operators $R_j$ and $R_j^T$, are denoted $\textbf{R}_j$ and $\textbf{R}_j^T$ respectively. This allows us to define the subdomain stiffness matrices restricted to $V_{h,0}( \Omega_j)$ as 
$\textbf{A}_j := \textbf{R}_j \textbf{A}\textbf{R}_j^T$ for $j = 1, \ldots, N$. In practice, we do not compute ${\bf A}_j$ from $\bf{A}$ via this double matrix product. Instead, we can equivalently assemble  ${\bf A}_j$ directly from the bilinear form \eqref{eqn:bilinear} on $\Omega_{j}$ with homogeneous Dirichlet boundary conditions on all artificial interior subdomain boundaries, i.e  all points ${\bf x} \in \partial \Omega_j$ that satisfy ${\bf x} \in \Omega_{j^{'}}$ for some other $j^{'} \neq j$.

\smallskip\noindent
The $1$-level Additive Schwarz method can then be defined as a preconditioner of \eqref{eqn:fe_matrix_system} via the operator
\begin{equation}\label{eqn:onelevelAS}
{\bf M}^{-1}_{AS,1} = \sum_{j=1}^N {\bf R}_j^T {\bf A}_j^{-1} {\bf R}_j
\end{equation}
Here, in the subscript, the $AS$ denotes additive Schwarz and the $1$ denotes a one-level method. In this case the preconditioner ${\bf M}^{-1}_{AS,1}$ approximates the inverse operator ${\bf A}^{-1}$ by a sum of local solves on overlapping subdomains, with homogeneous Dirichlet boundary conditions on interior boundaries. We will see in the numerical examples to follow that, for large problems, a single-level method is not sufficient, causing stagnation (high iteration counts) of the iterative solver. This stagnation of the iterative solver is caused by a few very small eigenvalues in the spectrum of the preconditioned problem. They are due to the lack of a global exchange of information in the preconditioner in the single-level method. A classical remedy is the introduction of a coarse grid problem that couples all subdomains at the second level \cite{Tos04}. To define our coarse problem, we introduce a coarse space $V_H \subset V_h$ (which we define below). We denote the restriction from the fine to the coarse space by the operator $R_H:V_h \rightarrow V_H$, with matrix representation ${\bf R}_H$. The two-level additive Schwarz preconditioner (in matrix form) is given by
\begin{equation}
{\bf M}^{-1}_{AS,2} = {\bf R}_H^T {\bf A}_H^{-1} {\bf R_H} + {\bf M}^{-1}_{AS,1} \quad \mbox{where} \quad {\bf A}_H = {\bf R_H} {\bf A} {\bf R_H}^T\,.
\end{equation}

\noindent\smallskip
Two natural questions arise: 
\begin{itemize}
\item What is a good choice of coarse space $V_H$ for composite applications? \vspace{-1ex}
\item How do we construct ${\bf A}_H$ efficiently on a distributed memory computer without assembling ${\bf A}$ directly?
\end{itemize}

\subsection{A Robust Coarse Space via Generalised Eigenproblems in the Overlaps (\textsc{GenEO})}

\noindent The ideal coarse space would capture the global low energy modes of ${\bf A}$ that jeopardise the performance of Krylov solvers. Specifically, in the two-level additive Schwarz setting, the modes not captured by the local solves are of interest. Yet, to compute those low-energy modes explicitly would be more expensive than inverting {\bf A} itself. Instead, the global low-energy modes can be approximated by stitching together local (optimal) approximations. These local approximations are solutions of specific {\sc Gen}eralised {\sc E}igenproblems in the {\sc O}verlaps, hence named \textsc{GenEO}, defined below. Importantly the local eigenproblems are independent and can trivially be computed in parallel. The robustness of \textsc{GenEO} has been proven for isotropic elasticity problems, Spillane et al. \cite{Spi14}, and numerically verified by the authors for anisotropic variants \cite{dunecomp}.

\smallskip\noindent The construction of the \textsc{GenEO} coarse space has two key steps: the definition of the generalised eigenproblems on the subdomains and the stitching together of the resulting local eigenmodes from each subdomain to form a global basis. This stitching process by means of {\em partition of unity} (PoU) operators is also incorporated in the local eigenproblems, therefore we construct the PoU operators first.

\begin{definition} (Subdomain Overlap). For each subdomain $\Omega_j$, the overlap region is defined by the set
$$
\Omega_j^\circ := \{ {\bf x} \in \Omega_{j}: \exists j^{\prime} \neq j \; \text{s.t.} \; {\bf x} \in \Omega_{j^{\prime}} \},
$$
i.e. the subset of $\Omega_j$ which belongs to at least one other subdomain.
\end{definition}

\begin{definition} \label{def:pou} (Partition of Unity).
The family of operators $\Xi_j : V_h(\Omega_j) \rightarrow V_{h,0}(\Omega_j)$, $j=1,\ldots,N$, defines a Partition of Unity if
\begin{eqnarray*}\label{eqn:pou}
  \sum_{j = 1}^N R_j^T \Xi_j ( v |_{\Omega_j}) = v, &  & \forall v
  \in V_h \,.
\end{eqnarray*}
\end{definition}
 
Since $R_{j'}^T \Xi_{j'} ( v |_{\Omega_{j'}}) = 0$ on  $\Omega_j \backslash \Omega_j^o$ for all $j' \not= j$, it follows from this definition that restricted to $\Omega_j \backslash \Omega_j^o$ each $\Xi_j$ has to be the identity operator.
In the overlaps, the choice of $\Xi_j$ is not unique. The simplest approach is to define $\Xi_j(v)$ such that each coefficient of the FE function $v$ is scaled by the number of subdomains the corresponding degree of freedom belongs to (see \cite{Spi14} for details). However, we also provide a smoother PoU as defined by Sarkis \cite{Sar97} in our implementation, but observe little difference in the performance of the overall solver (at most one iteration); we therefore keep the presentation here as simple as possible.

\smallskip
\noindent Given this set of local PoU operators $\Xi_j(\cdot)$, we can construct any global FE function $v_h \in V_h$ from local functions $v_h^{(j)} \in V_h(\Omega)$ as follows:
\begin{equation}
v_h= \sum_{j=1}^NR_j^T\Xi_j(v_h^{(j)})\,.
\end{equation}

\noindent In particular, we can define the local generalised eigenproblems that (once collected from each subdomain) provide the basis of the \textsc{GenEO} coarse space. The following definition can be rigorously motivated from theoretical considerations and we refer again to \cite{Spi14}. For each subdomain $\Omega_j$, $j = 1, \ldots, N$, we define the generalised eigenproblem: Find $(\lambda, p) \in \mathbb R^+ \times V_{h}(\Omega_j)$ such that
  \begin{equation}\label{eq-geneo}
   a_{\Omega_j} ( p, v) = \lambda a_{\Omega_j^o} ( \Xi_j ( p), \Xi_j ( v)),~~~~~ \forall v \in V_h ( \Omega_j),
  \end{equation}
where, for any $D \subset \Omega$, the bilinear form  $a_{D}$ is defined like $a$ in (\ref{eq-fe}) with the integral restricted to $D$.

\begin{definition}
(GenEO coarse space).
  For each subdomain $\Omega_k$ let $p_k^j$ be the
  eigenfunctions from  (\ref{eq-geneo}) with associated eigenvalues $\lambda_k^j$ in ascending order.
  Then, for some choice of $m_j \in \mathbb{N}$, the GenEO coarse space
  is defined as
  $$V_H := \mbox{span} \{ R_j^T \Xi_j ( p_k^j) ~:~ k = 1, \ldots, m_j,\, j = 1, \ldots, N \} . $$
\end{definition}

\noindent
The only parameter that remains to be chosen is the number of eigenmodes $m_j$ to be included in each subdomain $\Omega_j$. In order to ensure robustness and scalability of the solver, the condition number of the preconditioned system needs to be bounded from above, independent of $N$, $h$ and of the material properties. It has been shown in \cite{Spi14} that

\begin{equation}
 \kappa ( {\bf M^{- 1}_{{AS}, 2} A}) \leqslant C \max_{1 \leqslant j \leqslant N} \left( 1 +
     \frac{1}{\lambda_{m_j + 1}^j} \right) .
\end{equation}
where $C>0$ is a constant depending only on the geometry of the subdomains and where $\lambda_{m_j + 1}^j$ is the lowest eigenvalue whose eigenfunction is not added to the coarse space on $\Omega_j$. Thus, the desired robustness can be achieved by including all eigenfunctions in the coarse space whose eigenvalues are below an a priori chosen threshold. A particular threshold that turns out to provide an effective black-box choice for $m_j$ and also depends only on the geometry of the subdomain partition is to include all eigenfunctions with  $\lambda_k^j \le \text{diam}(\Omega_j) / \text{width}(\Omega_j^o)$.\footnote{For any $D \subset \Omega$, $\text{diam}(D)$ and $\text{width}(D)$ refer to the radius of the largest circumscribed and inscribed circle, respectively.} This simple threshold can be scaled by a constant factor, thus also scaling the condition bound of the preconditioned system by the same factor. As the number of iterations of the Krylov solver depends directly on the condition number, this allows us to balance the time spent in the iterative solver with the time spent on setting up the preconditioner.

\smallskip\noindent
The number of eigenfunctions that are used in the coarse space is problem-specific, but it turns out that for strongly structured coefficient distributions only a small number is typically sufficient. We will see in Section \ref{sec:examples} that the calculation of these local eigenmodes is not prohibitively expensive, while yielding excellent condition numbers and, due to the independence of the individual eigenproblems, parallel scalability.

\subsection{Implementation of \textsc{GenEO} on a High Performance Computer}

\noindent The two-level additive Schwarz preconditioner with \textsc{GenEO} coarse space is implemented within a collection of header files, which are located in \texttt{dune-pdelab} from \texttt{releases/2.6} and in the folder \texttt{solvers/geneo/} in prior releases. Here we describe our implementation. We are aware of only one other high performance implementation of \textsc{GenEO}, which can be found in the package \texttt{HPDDM} for which details are provided in {\em Jolivet et al} \cite{Jol12,Jol13}.

\smallskip\noindent
It is a main goal of such an implementation to fully exploit the excellent parallel scalability promised by the method's construction and theoretical properties. Therefore, each process $j$ will be assigned to subdomain $\Omega_j$ and only store relevant fine-level operators and functions in the form of local restrictions to $\Omega_j$. Further, per-subdomain stiffness matrix and eigenproblem solves will be run in parallel on the respective processes. Only scalable nearest-neighbour communication is needed, with the exception of setting up the coarse matrix, which consequently requires particular care.

\subsubsection {Partition of Unity (PoU) operator}
\noindent The partition of unity operator as defined by Def. \ref{def:pou} is stored locally on each processor (see Fig. \ref{fig:domaindecomp}). In practice, the partition of unity operator $\Xi_j$ is represented as a diagonal matrix $\textbf X^{(j)}$. In the simplest case, each diagonal entry of $\textbf X^{(j)}$ is set to one divided by the number of subdomains containing the associated degree of freedom, except for the subdomain boundary where entries are set to zero. Therefore, if ${\bf v}^{(j)}$ is a vector containing all nodal degrees of freedom of the FE function $v_h \in V_h(\Omega_j)$ in subdomain $\Omega_j$, the operation $\textbf X^{(j)} {\bf v}^{(j)}$ automatically maps $v_h$ into $V_{h,0}(\Omega_j)$. Such a PoU can be generated using existing parallel data structures in {\sc Dune} by adding a vector of ones 
and by enforcing both global and subdomain boundary conditions before and after communication. The implementation of the PoU operators is within the header file \texttt{geneo/partitionofunity.hh} under the function \texttt{standardPartitionofUnity(...)}. As the choice of partition of unity operator is not unique, we also provide the PoU in \cite{Sar97}, which is implemented in the same header file under the function \texttt{sarkisPartitionofUnity(...)}. This gives a `smoother' PoU operator, which is however restricted to {\color{black} equally distributed subdomain sizes}. Under testing, we noted no significant difference in the performance of the preconditioner when changing between the two different PoU operators.

\subsubsection{Subdomain eigenproblems}
\noindent The local generalised eigenvalue problems \eqref{eq-geneo} can be rewritten in matrix form as follows: Find eigenpairs $(\lambda^{(j)}_i,{\bf p}^{j}_i)$ with eigenvalues in ascending order and $\|{\bf p}^{j}_i\| = 1$ such that
\begin{equation}\label{eqm:matrixForm}
\textbf A_{\Omega_j} {\bf p}^{j}_i = \lambda^{j}_i \Big(\textbf X^{(j)} \textbf A_{\Omega^{\circ}_j} \textbf X^{(j)}\Big) {\bf p}^{j}_i\,, \quad \mbox{for} \quad j = 1, \ldots, N \quad \mbox{and} \quad i = 1, \ldots, m_j
\end{equation}
where $\textbf A_{\Omega_j}$ and $\textbf A_{\Omega^\circ_j}$ denote the stiffness matrices corresponding to the bilinear forms $a_{\Omega_j}(\cdot,\cdot)$ and $a_{\Omega_j^\circ}(\cdot,\cdot)$ on $V_h(\Omega_j)$ and $V_h(\Omega_j^o)$, respectively. They are solved using \texttt{ARPACK++} \cite{arpackpp}.

\smallskip \noindent
A customised wrapper has been developed to convert \textsc{Dune} data structures into a suitable format for \texttt{ARPACK++} \cite{arpackpp}. In order to regularise the problem, we employ \texttt{ARPACK++}'s {\em shift and invert spectral transformation mode} and, since we are interested in the smallest $m_j$ eigenvalues, we choose a small shift factor. The global coarse basis vectors $\boldsymbol{\Phi}_1, \ldots, \boldsymbol{\Phi}_{N_H} $ are obtained from the local eigenvectors by applying the PoU operator, i.e., $\boldsymbol{\Phi}_{i(j,k)} := {\bf X}^{(j)}{\bf p}^{j}_{k}$, and padding the rest of the global vector (outside $\Omega_j$) with zeros. Here, $(j,k) \mapsto i(j,k)$ is a one-to-one mapping from the local numbering of the eigenvectors on $\Omega_j$ to a global numbering, with $1 \le i(j,k) \le N_H = \sum_{j=1}^N m_j$.  

\subsubsection{Coarse space assembly }
\noindent The parallel assembly of the coarse system ${\bf A_H} = {\bf R}_H {\bf A} {\bf R_H}^T$ is not trivial in practice since process $j$ only has local access to rows and columns of ${\bf A}$ associated to degrees of freedom on sub-domain $\Omega_j$. We denote this submatrix ${\bf \tilde{A}}_j$. Note that ${\bf \tilde{A}}_j$ differs from the matrix ${\bf A}_j$ in \eqref{eqn:onelevelAS} in that it does not incorporate Dirichlet conditions on interior subdomain boundaries.

Furthermore the coarse space prolongation matrix ${\bf R}_H^T$ is only available in a distributed manner. Each basis vector $\boldsymbol{\Phi}_i$, $i \in \{1,\ldots,N_H\}$, is available only on process $j(i)$, where the unique $j(i) \in \{1,\ldots,N\}$ denotes the index of the subdomain $\Omega_{j(i)}$ associated with the eigenproblem \eqref{eqm:matrixForm} corresponding to $\boldsymbol{\Phi}_i$. However, due to the local support of the basis functions, one can break down the global matrix product into local products
\begin{equation}
({\bf A}_H)_{i,\ell} = ({\bf R}_H {\bf A} {\bf R}_H^T)_{i,\ell} = \left( \boldsymbol{\Phi}_i^T {\bf \tilde{A}}_{j(i)}  \right) \boldsymbol{\Phi}_\ell, \quad \mbox{for} \quad i,\ell = 1, \ldots, N_H,
\end{equation}
with a slight abuse of notation, denoting the local parts of the global vectors $\boldsymbol{\Phi}_i$ and $\boldsymbol{\Phi}_\ell$ restricted to $\Omega_{j(i)}$ again by $\boldsymbol{\Phi}_i$ and $\boldsymbol{\Phi}_\ell$. In the implementation, the matrix vector product in the bracket is local whereas the scalar product requires communication of (parts of)  vector $\boldsymbol{\Phi}_\ell$ from processor $j(\ell)$ to processor $j(i)$. This avoids having to communicate the local matrices ${\bf \tilde{A}}_{j}$.

\smallskip\noindent
We note that the locality of the basis functions implies $\boldsymbol{\Phi}_i {\bf \tilde{A}}_{j(\ell)} \boldsymbol{\Phi}_\ell = 0$ whenever $\Omega_{j(i)} \cap \Omega_{j(\ell)} = \emptyset$. Therefore, the parallel assembly of ${\bf A_H}$ requires only communication between processes assigned to overlapping subdomains. This locality of communication, which is demonstrated in Fig. \ref{fig:domaindecomp} (middle), can be exploited to set up ${\bf A_H}$ as a banded sparse matrix. Its parallel communication is 
implemented within \texttt{geneo/multicommdatahandle.hh}, allowing to pass basis functions between all processes at the same time and therefore make best use of available bandwidth. Combining sparsity and efficient communication, linear complexity in basis size can be achieved for this step.

\smallskip\noindent
Since the number of coarse degrees of freedom per process is too small, i.e., $m_j$ on process $j$, after assembly we distribute the resulting global coarse matrix ${\bf A}_H$ to all processors and duplicate the coarse solve on all processors to avoid fine-grained communication.  The communication of ${\bf A}_H$ is achieved directly with  \texttt{MPI} calls, as the \texttt{dune-pdelab} communication infrastructure is 
not designed for such global operations.

\smallskip\noindent
In case of the restriction of a distributed vector ${\bf v_h}$ representing a function $v_h \in V_h$, it follows that
\begin{equation}
({\bf R}_H{\bf v_h})_i = \boldsymbol{\Phi}_i^T {\bf v_h}.
\end{equation}
So, each row $i$ can be computed by the process associated with $\boldsymbol{\Phi}_i$, and the rows can be exchanged among all processes via \texttt{MPI\_Allgatherv}. Again, the communication effort increases with the dimension of $V_H$. On the other hand, the prolongation ${\bf R}_H^T {\bf v_H}$ of a vector ${\bf v_H}$ that is globally available on all processors, representing a $v_H \in V_H$, consists only of local contributions and hence can be computed in parallel without communication. 
{\color{black}
Since the result lies in $V_h(\Omega)$, the regular 
PDELab communication patterns can be used. This only involves communication between adjacent subdomains, making this a highly scalable process.}

\smallskip\noindent Each process executes its associated subdomain solve as well as the coarse space solve (redundantly). Where possible we use a sparse direct solver (\texttt{UMFPack}) \cite{Dav04}. For very large problems and a large number of parallel processors the coarse space becomes too large, and must itself be solved with a preconditioned iterative solver; in that case we use by default preconditioned CG with the \textsc{BoomerAMG} \cite{boomerAMG} preconditioner. It is important to note that in such cases since the coarse solve is inexact, the preconditioner for the (overall) Krylov method is now instationary. It is therefore necessary to switch from a standard preconditioned CG to a flexible Krylov solver. In our case we use Flexible GMRES as provided by \texttt{dune-istl} \cite{Bla07}.

\section{Using and Extending \texttt{dune-composites}} \label{sec:code}

\begin{figure}[ht!]
\includegraphics[width=0.6\linewidth]{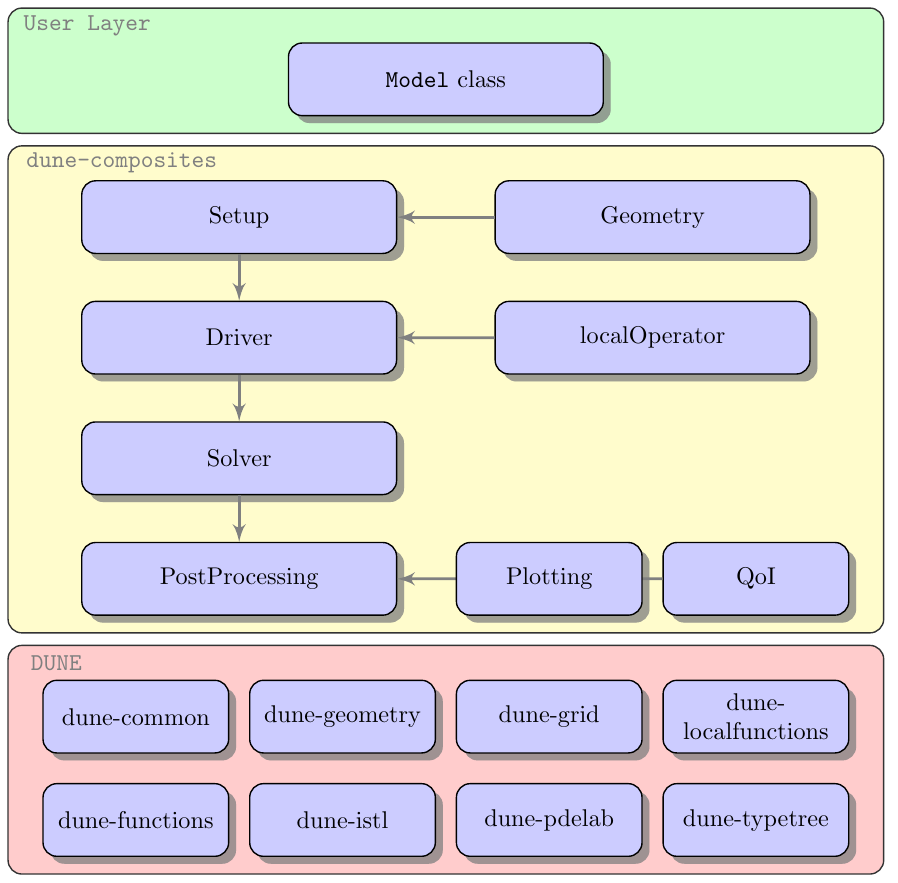}
\caption{\label{fig:structure}Code structure.}
\end{figure}

\noindent In this section we provide an overview of the \texttt{dune-composites} code, sufficient to enable other scientists to leverage the framework. The code is structured so that little additional knowledge of \textsc{Dune} and/or C++ is required to apply the code within the existing functionality. Figure \ref{fig:structure} shows the code structure. A user can extend any of the functionality e.g. implement a new solver, define a more complex nonlinear problem (e.g. cohesive zone) or introduce new types of elements. 

\smallskip\noindent
\subsection{Defining a \texttt{Model}}
At the highest level an analysis is defined by a user-defined \texttt{baseStructuredGridModel} class shown in green in Figure \ref{fig:structure}. 
This class defines all the key variables, functions and classes which describe the analysis, as well as storing any variables that are required for postprocessing or any later calculations. A base model class is provided, which can be inherited by each example. This provides default variables and functions, so that the user need only overwrite those functions which deviate from this base class. The \texttt{Model} class also defines the general loading on the structure and the boundary conditions, these include Dirichlet and Neumann conditions, but also thermal loading and multi-point constraints. Periodic boundary conditions are defined within the grid data structure, using \texttt{Model::LayerCake()}. Examples of user defined \texttt{Model} classes for a series of applications are provided in Sec. \ref{sec:examples} which follows.

\smallskip\noindent
\subsection{Internals of \texttt{dune-composites}}
 The functions provided in the subfolder \texttt{/Setup} provide support for the geometric setup of the grid geometry, material properties and boundary conditions. This includes the composite layering (or stacking sequence), the structural geometric shape of the component and in some cases adding a perturbation to the geometry to form a defect (see for example Sec. \ref{sec:Example1}). Because of uniform layering and planar anisotropy in composite laminates, our first version has focused on structured, overlapping grid implementations using \texttt{Dune::YaspGrid} and \texttt{Dune::GeometryGrid}. {\color{black} The \texttt{Dune::GeometryGrid} functionality allows us to apply any continuous transformation to the basic Cartesian mesh provided by \texttt{Dune::YaspGrid}.} Development of unstructured grid implementation using \texttt{Dune::UGGrid} \cite{UG} are described in the future work section \ref{sec:conclusions}.

\smallskip\noindent 
 The folder \texttt{/Driver} contains the key functions and classes which relate to the FE calculations beyond what is available directly from \texttt{Dune::PDELab}. In particular these include all  element calculations, the definitions of new FEs, solvers and preconditioners. The functions and classes are split between three folders:
\begin{itemize}
\item {\bf \texttt{/localOperators} } define the weak form of the equations to be solved on an element, along with any support functions. In our case for anisotropic linear elasticity equations we define the local operator \texttt{linearelasticity.hh} which returns the element stiffness matrix, load vector and residual as defined by equation \eqref{eq-fe}. 

\item {\bf \texttt{/FEM} } defines specialist finite elements beyond those defined by \texttt{Dune::PDELab}. In our case, these are the family of serendipity elements \cite{Arn11}. The use of these elements are then defined explicitly in the \texttt{Driver} class, where the FE space is set up on the grid.

\item {\bf \texttt{/Solvers}} defines specialist solvers and preconditioners beyond those defined in \texttt{Dune::PDELab::istl} \cite{Bla07}. The \texttt{Driver} uses the solver as defined by the \texttt{Model} class, defined by the templated class function \texttt{Model::solve()}. By default, as defined by \texttt{baseStructuredGridModel}, for parallel calculations we use a CG Krylov solver, preconditioned with either a one or two level additive Schwarz method. Two-level methods use \textsc{GenEO} as the coarse space as long as \texttt{ARPACK++} is available. If not, the coarse space consists of only the zero energy modes \cite{Tos04}. In sequential mode, in particular for the coarse and the local solves, a sparse direct solver \texttt{UMFPack} \cite{Dav04} and the iterative CG preconditioned with AMG as provided by \texttt{Dune::PDELab::ISTLBackend\_SEQ\_CG\_AMG\_SSOR}, are also available. In \texttt{/Solvers/hypre} we provide a wrapper to the external parallel solvers provided by \texttt{hypre} \cite{Fal06}, including \texttt{boomerAMG} \cite{boomerAMG}.
\end{itemize}
The analysis is defined via a \texttt{Driver}. This is a class which provides the complete solution procedure, from the setup of the grid, the definition of a finite element space, assembly of the global stiffness matrix and load vector, and finally, calling the solver, followed by the Post Processing routines, as defined by the \texttt{Model} class. In this version, we provide two driver classes for the examples considered in Sec. \ref{sec:examples}: \texttt{/FEMDriver/linearStatic.hh} and \texttt{/FEMDriver/ThermalStatic.hh}.  

\smallskip\noindent
The folder \texttt{/PostProcessing} contains various classes and functions for the postprocessing of results. By default, after the solution has been calculated the driver initiates certain post processing steps, in particular the calculation of stresses, the creation of the necessary data for any plots and finally the calculation of any further quantities of interest. All of these routines can be modified by the user. 

\section{\texttt{dune-composites} - Examples}\label{sec:examples}

\noindent In this section, we introduce and demonstrate the functionalities of \texttt{dune-composites} using a series of examples of increasing complexity. The examples are intended as a starting point for researchers implementing their own studies, whilst also demonstrating the significant computational gains \texttt{dune-composites} is able to achieve in comparison to the commercial package \textsc{Abaqus} \cite{Abaqus}.

\smallskip\noindent
To simplify the definition of the examples in the following we assume that all cases use the same material properties. The orthotropic fibrous layers are assumed to be of thickness $t_p = 0.23$mm, with elastic moduli 
\begin{align}
E_{11} = 162 \mbox{GPa}, \quad E_{22} = E_{33} =  10 \mbox{GPa},& \quad G_{12} = G_{13} = 5.2\mbox{GPa}, \quad
 G_{23} =  3.5\mbox{GPa}, \\ \nonumber
 \nu_{12} = \nu_{13} = 0.35\quad &\mbox{and} \quad \nu_{23} = 0.5, \nonumber
\end{align} 
whereas the isotropic resin rich layers are assumed to be $t_i = 0.02$mm thick, with isotropic properties $E = 10$GPa and $\nu = 0.35$. These particular values are taken from a previous study by the authors \cite{Fle16}.

\subsection{Example 1: A Flat Composite Plate}\label{sec:Example1}

\noindent For the first two examples we consider a flat composite plate $[0,100\mbox{mm}]\times[0,20\mbox{mm}]$ under various loading conditions. 
The laminate is made up of $12$ identical composite layers arranged in the following composite stacking sequence 
\begin{equation}\label{eqn:stackingSeq_Example01}
  \left[\mp 45^\circ/0^\circ/90^\circ/\pm 45^\circ/\mp 45^\circ/90^\circ/0^\circ/\pm45^\circ\right].
\end{equation}
The composite layers are separated by $11$ isotropic resin interface layers, giving a total thickness of $T = 2.98$mm. In each of the examples, we discretise the geometry with quadratic, 20-node serendipity elements (with full Gaussian integration). For the base mesh, which will be refined, we take $20$ elements in the $x$-direction, $5$ in the $y$-direction and through thickness $2$ per composite layer and $1$ per interface layer. This gives a total number of $3,500$ elements, with $13,608$ degrees of freedom.

\smallskip\noindent
The geometry, the stacking sequence and the initial finite element mesh are introduced into a model by overwriting the base class function \texttt{Model::LayerCake()}, with the following user defined function.
\begin{lstlisting}
void inline LayerCake(){
   std::string example01a_Geometry = "stackingSequences/example1a.csv";
   LayerCakeFromFile(example01a_Geometry);
   GeometryBuilder();
}
\end{lstlisting}
Here the geometry and grid are defined by a file \texttt{"stackingSequences/example1.csv"}.

\smallskip
\noindent In these first two examples, we demonstrate a very simple setup and run on a single processor as well as on a few processors. We consider a cantilever beam with a uniform pressure of $0.01$ MPa applied to the top face and the following boundary conditions:
\begin{equation}
\label{eq:bc_ex1}
u_1 = u_2 = u_3 = 0 \quad \mbox{at} \quad x = 0 \quad \mbox{and} \quad \sigma_{33} \cdot {\bf n}_3 = -q \quad \mbox{at} \quad z = T.
\end{equation}
All other boundary conditions are assumed to be homogeneous Neumann conditions, i.e 
\begin{equation}
\sigma_{ij}\cdot n_j = 0.
\end{equation}
Boundary conditions are implemented by overwriting the two class functions \texttt{Model::isDirichlet()} and \texttt{Model::evaluateNeumann()} as follows
\begin{lstlisting}
bool inline isDirichlet(FieldVec& x, const int i){
   return (x[0] < 1e-6);
} 
\end{lstlisting}
{\color{black} Here $i$ refers to the direction being restricted, see \eqref{eq-boundary}.}
\begin{lstlisting}
inline void evaluateNeumann(const FieldVec& x, FieldVec& h,
                            const FieldVec& normal) const{
    h = 0.0; // initialise to zero
    double T = R[0].L[2]; // Thickness
    if (x[2] > T - 1e-6){
        h[2] += q;
    }
}
\end{lstlisting}
We note that \texttt{Model::evaluateDirichlet()} need not be overwritten since by default it returns homogeneous boundary conditions (i.e. ${\bf u}({\bf x}) = {\bf 0}$ ) for all those points $x$ marked as Dirichlet boundary conditions by \texttt{isDirichlet()}. Furthermore, by default loading under the weight of the structure is included by providing density as an input parameter. We do not wish to include it in this example and therefore we must also overwrite the function \texttt{Model::evaluateWeight()}
\begin{lstlisting}
inline void evaluateWeight(FieldVec& f, int id) const{
    f = 0;
}
\end{lstlisting}
{\color{black} Here $f$ is the output density in a given element $id$.}

\subsubsection{Example 1a: A flat composite plate -- getting started}
\noindent
Our first study computes the maximum vertical deflection of the cantilever beam as our quantity of interest $$Q({\bf u}) = \max_{{\bf x} \in \Omega} u_3({\bf x}).$$
This is done by providing the following user defined function
\begin{lstlisting}
template<class GO, class V, class GFS, class C>
void inline postprocess(const GO& go, V& u, const GFS& gfs, const C& cg){
    using Dune::PDELab::Backend::native;
    double local_u3_max = 0.0;
    for (int i = 0; i < native(u).size(); i++){ // Loop over each vertex
        double u3 = std::abs(native(u)[i][2]);
        if (local_u3_max < u3) { local_u3_max = u3; }
    }
    MPI_Allreduce(&local_u3_max, &QoI, 1,
                  MPI_DOUBLE, MPI_MAX, MPI_COMM_WORLD);
}
\end{lstlisting}
This function loops over the solution at each vertex \texttt{native(u)[i]}, and records the maximum vertical displacement \texttt{native(u)[i][2]}. Since for a parallel run, this maximum is only the maximum on the local subdomain associated with a given processor, the final command \texttt{MPI\_Allreduce()} finds the maximum vertical displacement over all subdomains (processors). The final result is stored in \texttt{QoI}, {\color{black} a member of the \texttt{baseStructuredGridModel} class}.

\smallskip\noindent
We use the default sequential and parallel solvers. On a single processor (e.g.~with the call \texttt{./Example1a}) the sparse direct solver \texttt{UMFPack} \cite{Dav04} is used. Otherwise, if more than one processor is used (e.g.~with the call \texttt{mpirun -np 8 ./Example1a}) the equations will be solved with CG, preconditioned with a one-level additive Schwarz preconditioner, as defined by \eqref{eqn:onelevelAS} using \texttt{UMFPack} as the local solver on each subdomain.

\begin{figure}
\centering
\includegraphics[width = 0.8\linewidth]{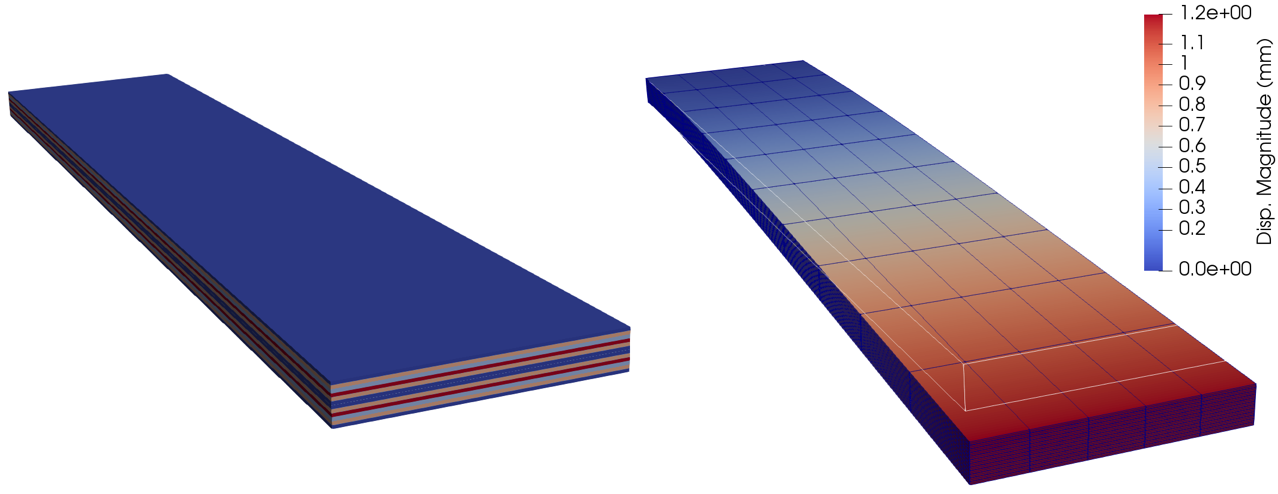}
\caption{\label{fig:Example01a}Visualisation of results for \texttt{Example01a} using \textsc{Paraview} (left) Visual output of laminate and stacking sequence using \texttt{plotProperties()} function (right) Visualisation of solution, in deformed coordinates (scalar factor of displacement is $4$). }
\end{figure}

\smallskip\noindent As output, the quantity of interest is printed to the screen (\texttt{Maximum Vertical Displacement in\linebreak Example01a = 1.23992mm}). Furthermore, the data for plots of the laminate stacking sequence, the solution (deformation) and the stress field are generated and  provided in three files named \texttt{Example01a\_xxx.vtu}. The stacking sequence and solution are shown in Fig. \ref{fig:Example01a}.

\subsubsection{Example 1b: A flat composite plate -- testing preconditioners (up to 32 cores)}

\noindent 
In \texttt{Example01b}, we test our new preconditioner \textsc{GenEO} on up to 32 processors. In this example we also demonstrate the inclusion of a failure criterion. To do this we change the quantity of interest to be the pressure $q = q^\star$ in the boundary condition \eqref{eq:bc_ex1} at which the laminate fails according to the Camanho criterion \cite{Cam03}, defined by the functional
\begin{equation}\label{eqn:Hashin}
\mathcal F(\sigma ({\bf x})) = \sqrt{\left(\frac{\sigma^+_{33}}{s_{33}}\right)^2 + \left(\frac{\sigma_{13}}{s_{13}}\right)^2 + \left(\frac{\sigma_{23}}{s_{23}}\right)^2 }.
\end{equation}
We apply the Camanho criterion only in the resin-rich interface layers and we say that failure occurs at a  load $q^\star$ if $\max_{{\bf x} \in \Omega^{\text{Inter}}} \mathcal F(\sigma ({\bf x})) = 1$. However, since the problem is linear it suffices to solve only one problem with an arbitrary load $q$. The failure load is then given by $q^\star := q / \max_{{\bf x} \in \Omega^{\text{Inter}}} \mathcal F(\sigma({\bf x))}$. Expression \eqref{eqn:Hashin} is implemented in the file \texttt{PostProcessing/FailureCriterion/Camanho.hh} within the code. 
Within \texttt{linearStaticDriver}, by default, the stress field (per element) is stored within the container \texttt{stress\_mech} (a $6 \times 1$ vector). To compute $q^\star$, the Camanho functional is first calculated in each element. The maximum is then found by once again overwriting the class function \texttt{Model::postprocess}. The material allowables, $s_{33} = 61$ MPa, $S_{13} = 97$ MPa and $s_{23} = 94$ MPa, in \cref{eqn:Hashin} are stored in a \texttt{std::vector<double> p}.

\hspace*{-\parindent}%
\begin{minipage}{\linewidth}
\begin{lstlisting}
template<....>
void inline postprocess(...){
    //material allowables in MPa
    const std::vector<double> p = {61., 97., 94.};
    double Fm = 0.0;
    for (int i = 0; i < stress_mech.size(); i++){
        double F = Camanho(stress_mech[i], elemIndx2PG[i], p);
        if (Fm < F) { Fm = F; }
    }
    double Fm_all;
    MPI_Allreduce(&Fm, &Fm_all, 1, MPI_DOUBLE, MPI_MAX, MPI_COMM_WORLD);
    Q = pressure / Fm_all; // Failure load
}
\end{lstlisting}
\end{minipage}

\smallskip\noindent
Different failure criteria can be implemented by defining other user-defined functionals of the stress tensor, similar to \texttt{Dune::Composites::Camanho()}. In this simple test, we note that failure initiates due to high through thickness stresses in the interface between layers ($\sigma_{13}$ and $\sigma_{23}$) as the laminate bends. For further engineering discussion of the failure of composites under the Camanho criterion we point the reader to the original paper \cite{Cam03} and to \cite{Fle16,dunecomp}.

\begin{figure}
\centering
\includegraphics[width=0.58\linewidth]{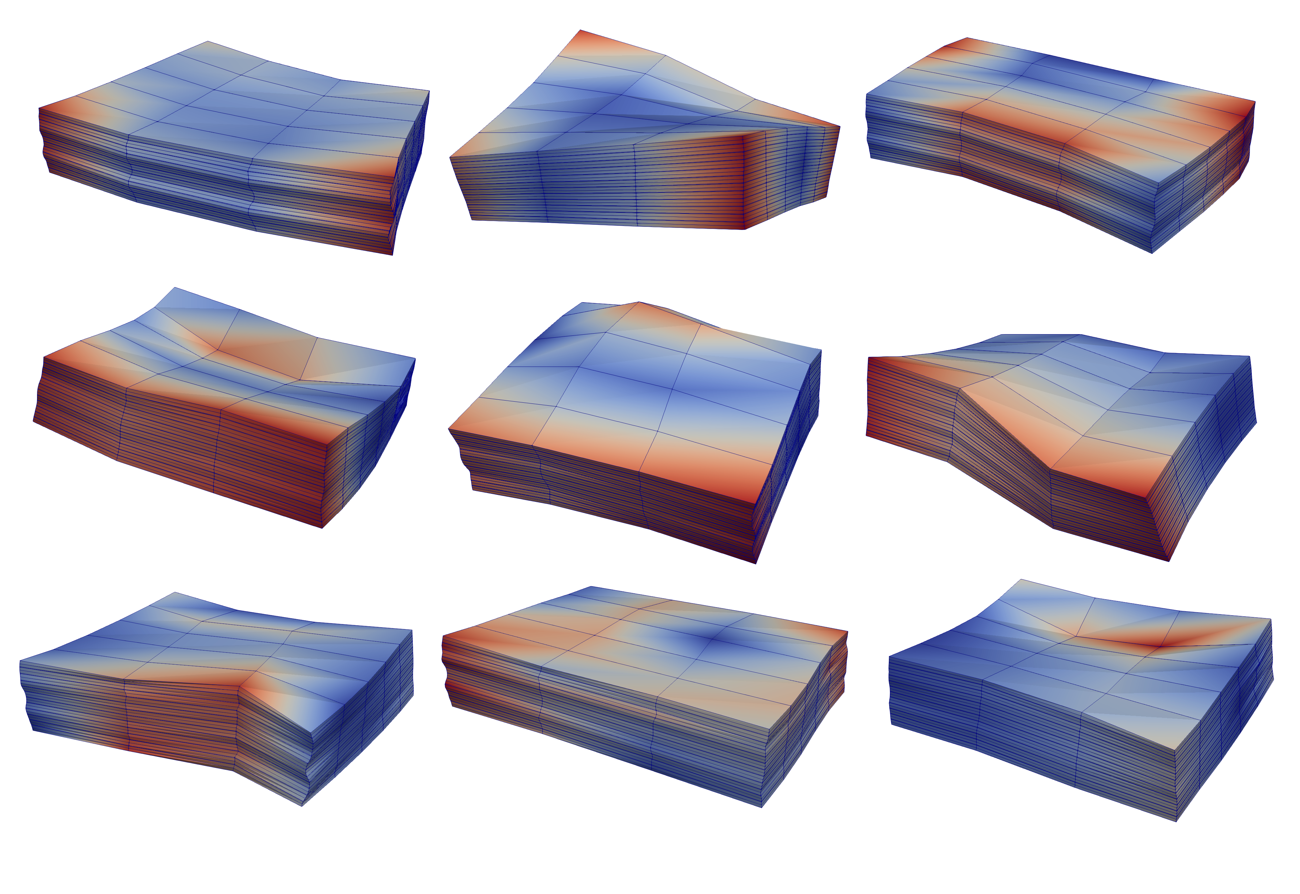}
\includegraphics[width=0.4\linewidth]{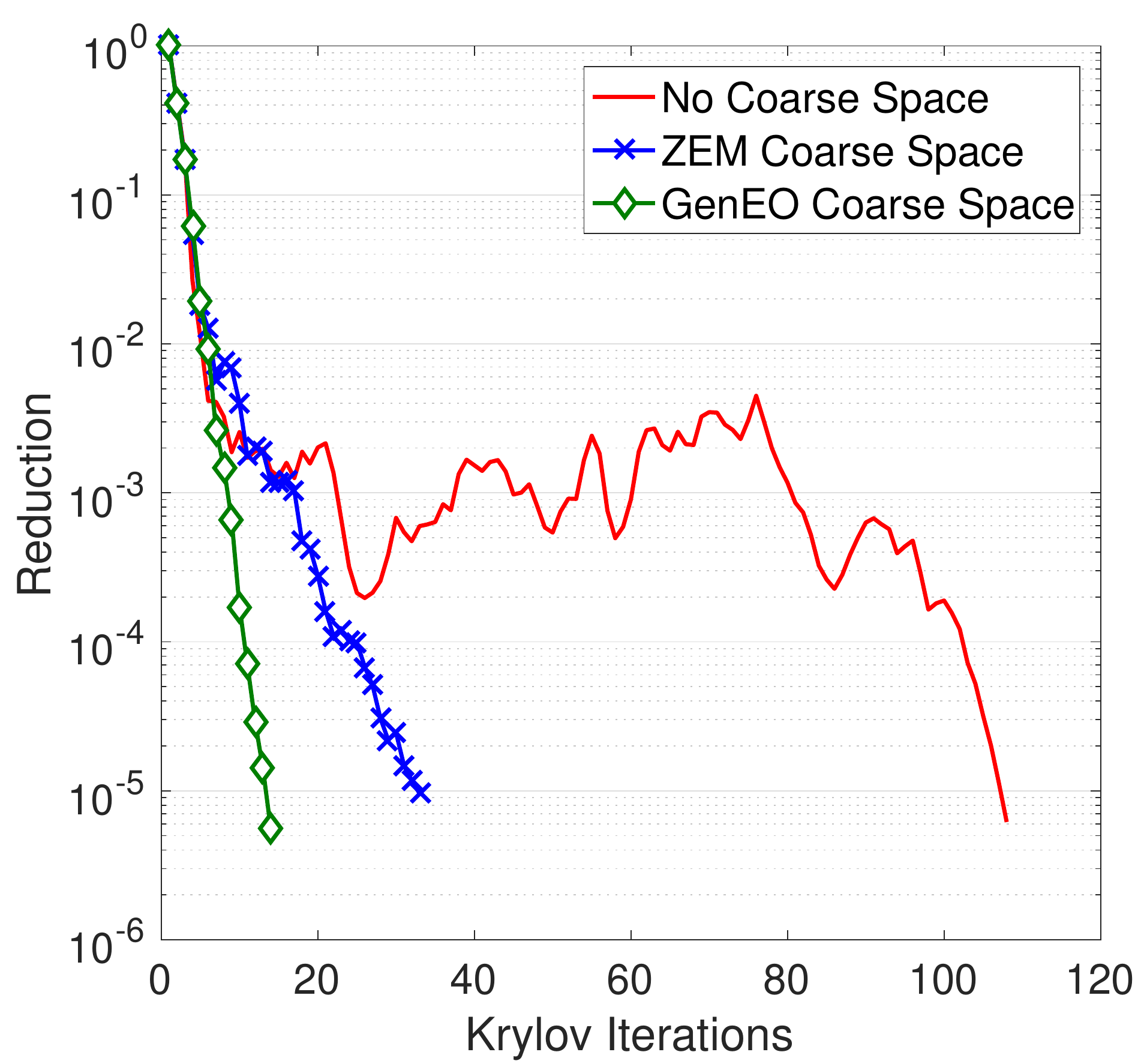}
\caption{\label{fig:iterations_comparision}(Left) The eigenvectors corresponding to the first nine non-zero eigenvalues on a subdomain with no global Dirichlet boundary. (Right) The reduction of the residual against CG iterations for \texttt{Example01b} using no coarse space, only zero energy modes (ZEM) and the full \textsc{GenEO} coarse space.} 
\end{figure}

\smallskip\noindent
We use this test example to demonstrate the influence of the GenEO coarse space on the parallel iterative solver. For the first experiment we use 16 processors. Fig.~\ref{fig:iterations_comparision} (left) shows the first nine non-zero energy modes of a subdomain with no global Dirichlet boundary. Linear combinations of these functions together with the zero energy modes (six rigid body translations and rotations) provide a good low dimensional representation of the system on that subdomain consisting of those modes most easily energetically excited. Fig.~\ref{fig:iterations_comparision} (right) shows the influence of the coarse space on the number of iterations for the preconditioned Krylov Solver (pCG), comparing no coarse space (one-level additive Schwarz), only the zero energy modes (ZEM) and the GenEO coarse space. The need for a coarse space is clear; with no coarse space, even in this simple test case we observe the well-documented stagnation phenomenon for iterative solvers \cite{Tos04} between Iteration $10 - 100$. With a coarse space (ZEM or GenEO), the convergence shows no stagnation and it is much faster -- close to optimal with \textsc{GenEO}.

\smallskip\noindent
Next we want to study the robustness of \textsc{GenEO} as a function of the number of subdomains in comparison to one-level additive Schwarz (AS) and ZEM. We consider a fixed size problem and increase the number of subdomains. We note that the tests can be run with
\begin{center}
\texttt{mpirun -np 16 ./Example01b} $\quad$ or $\quad$ \texttt{mpirun -np 16 ./Example01bBoomerAMG}
\end{center}
In each case we record the condition number, the dimension of the coarse space $\mbox{dim}(V_H)$ (if applicable) and the number of CG iterations to achieve a residual reduction of $10^{-5}$. The results are summarised in Table~\ref{tab:Example01b}. We see that the iteration counts (and the condition number estimates) increase steadily with the number of subdomains when no coarse space is used. The condition number estimate is still fairly big if only the zero energy modes are used and the iteration counts also increase steadily with the number of subdomains. In contrast, the iterations and the condition number estimates remain constant for the \textsc{GenEO} preconditioner. We also add a comparison with \texttt{boomerAMG} \cite{boomerAMG} for this test problem.
\texttt{BoomerAMG} provides a large number of parameters to fine-tune. We retained the defaults for most parameters (HMIS coarsening without aggressive refinement levels and a hybrid Gauss-Seidel smoother). We used blocked aggregation with block size $3$ as recommended for elasticity problems. A strong threshold of $0.75$ was chosen after testing values in the range from $0.4$ to $0.9$.
Due to a lower setup cost with this parameter setting, the \texttt{boomerAMG} solver is faster in actual CPU time, but the numbers of iterations -- albeit also constant -- are more than $10\times$ bigger. For more complex geometries, \texttt{boomerAMG} does not perform very well and in our tests it does not scale beyond about $100$ cores in composite applications.

\begin{table}[t]
\centering
\begin{tabular}{|c|c|c|c|c|c|c|c|c|c|c|}
\hline
\multirow{2}{*}{$N$} & \multicolumn{2}{c|}{AS} & \multicolumn{3}{c|}{ZEM} & \multicolumn{3}{c|}{GenEO} &  \multicolumn{2}{c|}{\textsc{BoomerAMG}} \\ \cline{2-11} 
                            & it       & cond $\kappa$         & it   & cond $\kappa$     & dim($V_H$)   & it   & cond $\kappa$  & dim($V_H$)      & it       & Num. levels \\ \hline
4                           & 89      & 79,735       & 26   & 394      & 12           & 16   & 10    & 78              & 258     & 10     \\
8                           & 97      & 84,023       & 30    & 245        & 42           & 15   & 9     & 126             & 258     & 11     \\
$16^\star$                   & 107      & 98,579       & 36   & 177      & 84           & 16   & 10    & 182             & 257     & 12     \\
32                          & 158      & 226,871      & 42   &  230     & 168          & 16   & 9     & 526             &  263    &   12    \\ \hline
\end{tabular}
\caption{\label{tab:Example01b}Demonstration of performance of different preconditioners for \texttt{Example01b} for fixed problem size (30,000 DOFs) but increasing the number of subdomains: Number of pCG iterations (it), coarse space
dimension (dim$(V_h)$), an estimate of the condition number $\kappa$ {\color{black} of the preconditioned system matrix $\bf M_{AS,2}^{-1} A$}.}
\end{table}

\subsection{Example 2 : Corner Unfolding -- Validation \& Performance Comparison with \textsc{Abaqus} (up to 32 cores) }

\noindent This example is motivated by the industrial challenge of certifying the corner-bend strength of a wingspar as its corner unfolds due to the internal fuel pressure in an airplane wing. We use this example to demonstrate the validity of the results of \texttt{dune-composites} by comparing the stresses computed with those given by \textsc{Abaqus}. We also make a cost comparison between the two software packages up to $32$ cores.

\begin{figure}
    \centering
    \includegraphics[width=\linewidth]{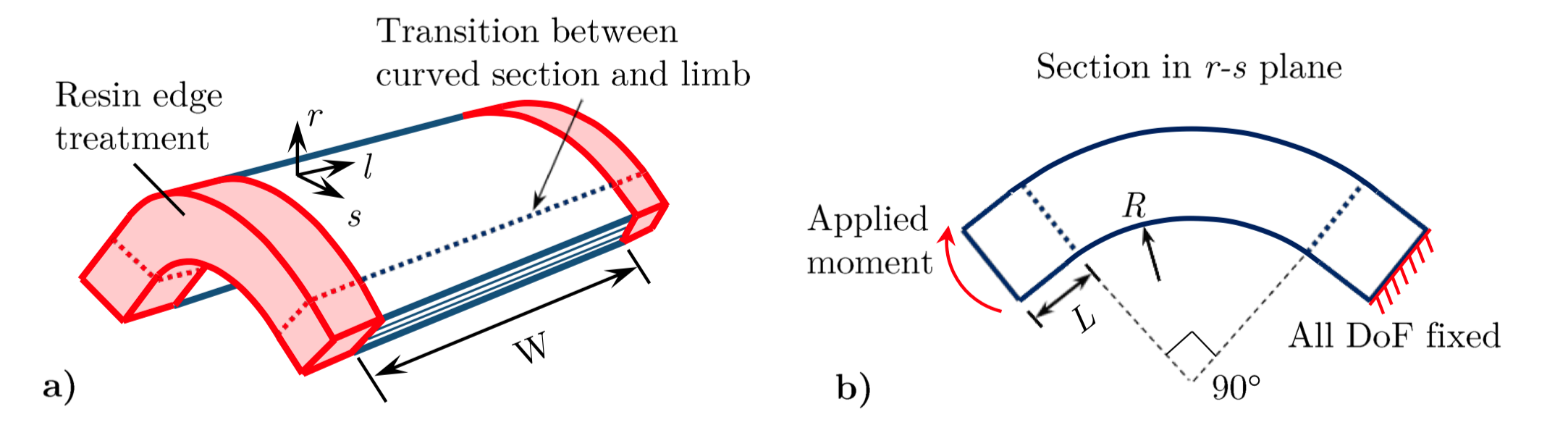}
    \caption{\label{fig:dune-vs-abaqus-setup}(Left) Diagram of the corner bend specimen with resin edge treatment. (Right) Cross section of the corner showing the loading conditions.}
  \end{figure}
 
\smallskip\noindent
The model setup is shown in Fig. \ref{fig:dune-vs-abaqus-setup}. We consider the curvilinear coordinate system ($s,r,\ell$), where $s$ is around the radius, $r$ is outwards (or normal) to the laminate and $\ell$ runs along the length of the sample. For our particular test, the two limbs of the coupon are of length $L = 3$mm and border a corner with a radius of $R = 6.6$mm. The width is taken to be $W = 15$mm. The 12 plies and the 11 interfaces have the same properties as in \texttt{Example01}, but the stacking sequence is slightly different, given by
\begin{equation}\label{eqn:stackingSeq_Example02}
  \left[\mp 45^\circ/90^\circ/0^\circ/\mp 45^\circ/\mp 45^\circ/0^\circ/90^\circ/\pm45^\circ\right].
\end{equation}
 Furthermore, we apply a resin treatment of $2$mm to the free edges of the laminate as shown in Fig. \ref{fig:dune-vs-abaqus-setup} (left). The advantages of edge treatment have been shown in \cite{Fle16}. It reduces conservatism in the design of aircraft structures, as well as making the analyses more reliable by eliminating stress singularities at the free edges.

\smallskip\noindent
Away from the points of contact, a standard four-point bend test as detailed in ASTM D6415 \cite{ASTM} generates a pure moment on the corner. To simulate such a moment, all degrees of freedom at the boundary of one limb are clamped. At the other end, all nodes are tied with a multipoint constraint where a running moment of $96.8$ Nmm/mm is applied. This is achieved by applying an offset load from the mid-plane of the laminate as a Neumann boundary condition, implemented with the user-defined function \texttt{evaluateNeumann()}.

\smallskip\noindent
The finite element mesh consists of $56 \times 56$ columns of hexahedral 20-node serendipity elements (element \texttt{C3D20R} in ABAQUS, \cite{Abaqus}) in its local $l$ and $s$ coordinates. In the $r$ direction, each (fibrous and resin) layer is discretised into 6 elements, leading to an overall number of $432,768$ elements. All of the geometry and mesh parameters are defined in \texttt{stackingSequenes/example2.csv}. To ensure a sufficient resolution of the strong gradients of the solution at the free edges and at the material discontinuities, the mesh is graded along the width towards both free edges and in the radial direction towards each of the fibre-resin interfaces. Grading is defined by the ratio between largest and smallest elements in the mesh, called the bias ratio and chosen to be 400 between the center and the edge in the $l$ direction and 10 between the layer centers and interfaces in the $r$ direction. The specification of the geometry and of the mesh grading can be defined in the function \texttt{gridTransformation()}. 

\begin{figure}
\includegraphics[width = 0.59\linewidth]{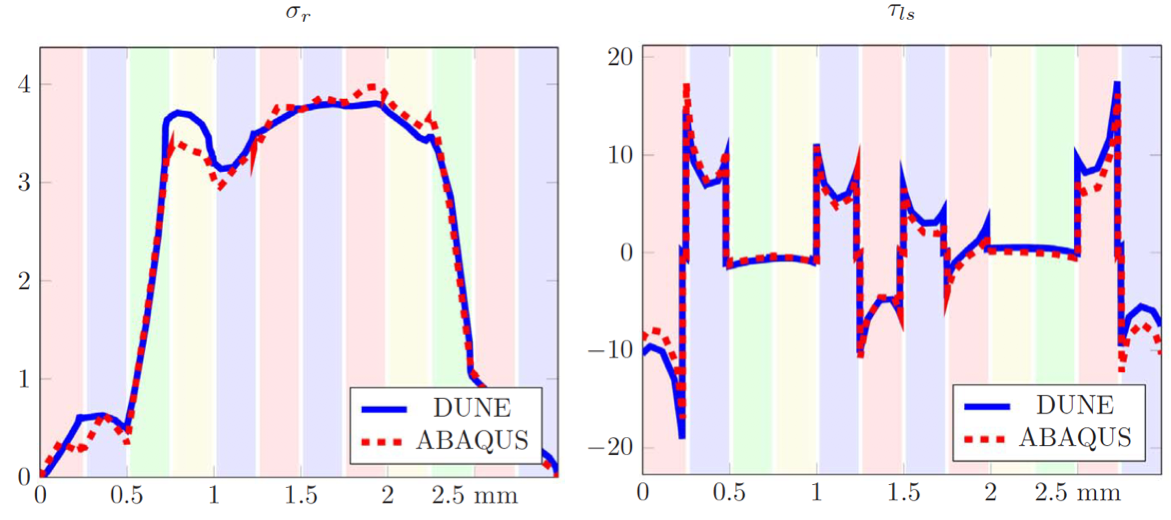}
\includegraphics[width = 0.32\linewidth]{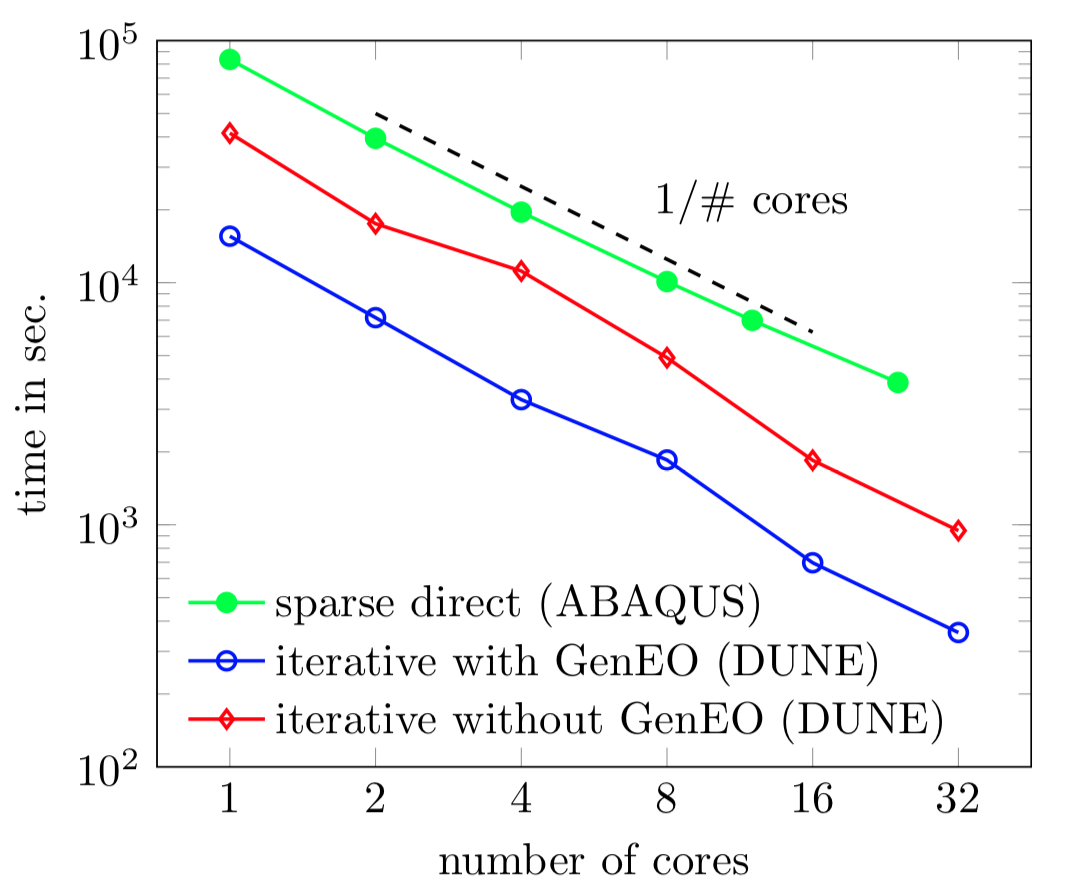}
\caption{\label{fig:duneabaqusstress}(Left \& Middle) Stresses (in MPa) as functions of the distance r 
from the outer radius at the apex of the curve, at 2.156mm from the edge
of the resin-edge-treated laminate (\texttt{dune-composites}, solid blue;
\textsc{Abaqus}, dotted red). The background colours indicate the
stacking sequence: $+45^\circ$=red, $-45^\circ$=blue, $90^\circ$=green,
$0^\circ$=yellow. (Right) Cost comparison between the sparse direct solver implemented in \textsc{Abaqus} and the iterative preconditioned CG solver in \texttt{dune-composites}.}
\end{figure}

  \smallskip\noindent
In Fig.~\ref{fig:duneabaqusstress} (left \& middle), we compare the radius stress (denoted by $\sigma_r$) and the through-thickness shear stress (denoted by $\tau_{s\ell}$) recovered from \texttt{dune-composites} and \textsc{Abaqus}. We see good agreement between the two codes. There are two small differences: \textsc{Abaqus} uses reduced integration while our example uses full integration and the stresses are not recovered in an identical way from the displacements in the two codes. In Fig.~\ref{fig:duneabaqusstress} (right) we see the absolute cost and the parallel scalability of the sparse direct solver in \textsc{Abaqus} and the iterative CG solver in \texttt{dune-composites} for a fixed total problem size, i.e. a strong scaling test. The red and blue curves show one-level and two level overlapping Schwarz methods respectively, both of which perform better than the sparse direct solver (green) in \textsc{Abaqus}. However, both codes show optimal parallel scalability up to 32 cores. 
In \texttt{dune-composites} the problem is decomposed into 8 subdomains for $1-8$ cores, distributed evenly to the available cores. On 16 and 32 cores, each core is passed exactly one subdomain, i.e., the number of subdomains is 16 and 32, respectively. The local problems on each subdomain are solved using the sparse direct solver \texttt{UMFPack} \cite{Dav04}. The simulations in \textsc{Abaqus} are with a parallel sparse direct solver, based on a parallel multi-frontal method similar to \cite{mumps}. \textsc{Abaqus}'s iterative solver, which is based on CG  preconditioned with \texttt{ML} \cite{ML}  (another black-box AMG preconditioner), does not converge in a reasonable time for this problem.  Therefore the computational gains observed here are really the difference between using a direct and a robust iterative solver. Importantly, we note that the parallel sparse direct solver, available in \textsc{Abaqus} does not scale beyond 64 cores \cite{Abaqus}, making it unsuitable for problems much bigger than that considered here, and reinforcing the need for robust iterative solvers and therefore \texttt{dune-composites} as a package.

\subsection{Example 3 : Large Composite Structure -- Parallel Efficiency of \texttt{dune-composites} (up to 15,360 cores)}\label{sec:wingbox}

\smallskip \noindent 
 The industrially motivated problem described in this section is to assess the strength of a wingbox with a small localised wrinkle defect. Wrinkle defects, which can form during the manufacturing process \cite{Dod14,Bel17}, occur at the layer scale. They lead to strong local stress concentrations \cite{dunecomp, San18}, causing premature failure. Naturally, good mesh resolution around the defect is required, leading to finite element calculations with very large number of degrees of freedom. We leave the engineering discussion of the results to a future engineering publication, using it instead to demonstrate both weak and strong scalability of \texttt{dune-composites} up to 15,360 subdomains. The experiments in this section were performed using the UK national HPC cluster \textsc{Archer},  which has $4,920$ Cray XC30 nodes with two 2.7 GHz, 12-core E5-2697 v2 CPUs each.

\begin{figure}[t]
\centering
\includegraphics[height = 3 cm]{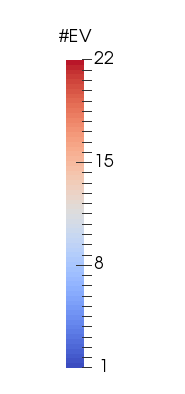}\hspace{-1.8em}
\includegraphics[height = 4.cm]{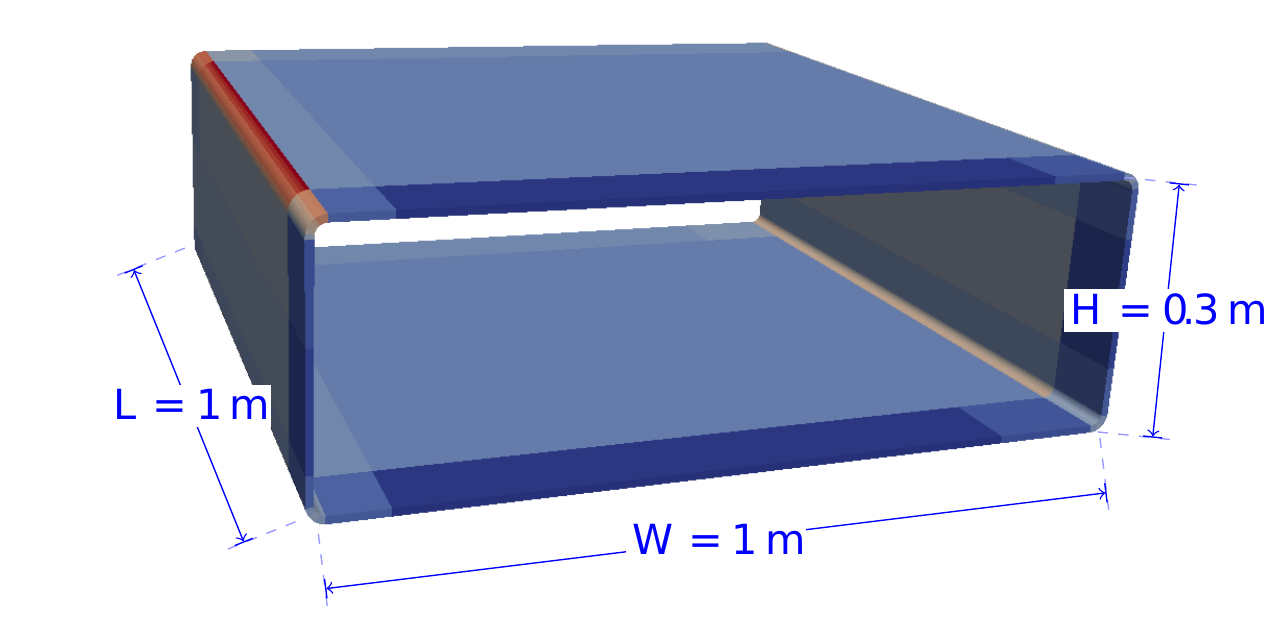}\hspace{1em}
\includegraphics[height = 4.cm]{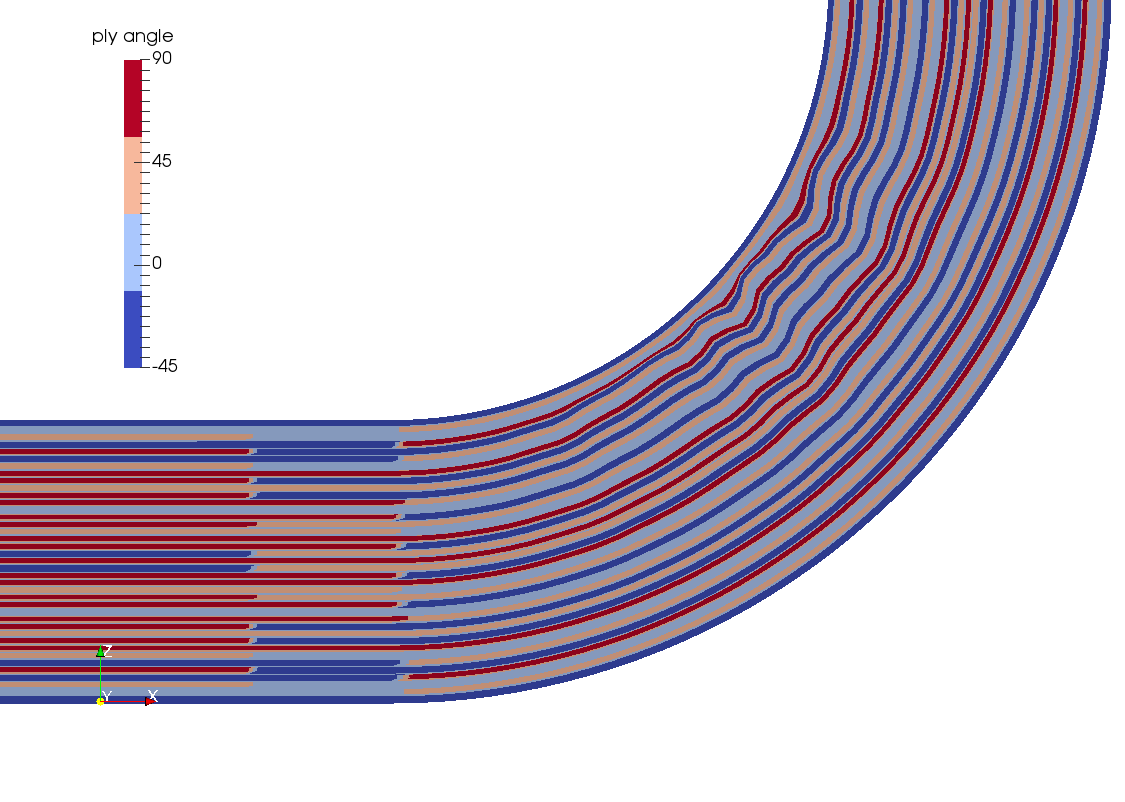}
\caption{\label{fig:wingboxgeom}(Left) Geometry of the wingbox with dimensions; the colouring shows the number of eigenmodes used in \textsc{GenEO} in each of the subdomains of Setup 6 in Table~\ref{table-wingbox}. (Right) Close-up plot of the corner of the wingbox using \texttt{plotProperties()}, which shows the wrinkle and the inter-lacing of the different stacking sequences in the corner, cover and spar regions.}
\end{figure}

\smallskip\noindent For these tests we model a single bay of a wingbox of width $W = 1$m, height $H = 300$cm and length $L = 1$m, as shown by the schematics in Fig.~\ref{fig:wingboxgeom} (left). The laminates were assumed to be of constant thickness, made up of 39 composite layers (as well as 38 interfaces) giving a total thickness of $T=9.93$cm with an internal radius of $15$mm in the corners. As in a typical aerospace application, the stacking sequence differs in the covers (top and bottom), corners and in the spar (sides) with the following approximate percentage breakdowns of $0^\circ$, $\pm45^\circ$ and $90^\circ$:
\begin{equation}
[50\%, 40\%, 10\%] \quad \mbox{(covers);}  \quad [20\%, 60\%, 20\%] \quad \mbox{(corners);}  \quad \mbox{and} \quad [15\%, 70\%, 15\%] \quad \mbox{(spars)}.
\end{equation}
We reiterate that this example serves to represent structural scale modelling. Therefore, sub-structural phenomena, such as stiffening of the upper and lower covers, are not modelled here. The specific layer-sequencing has been chosen, using a discrete optimiser, to ensure that each laminate is balanced, symmetric with no bend-twist coupling, whilst maximising the number of continuous orientations around the wing box. Transitions between each of the stacking sequences are achieved over a relatively short segment of 5cm, and the chosen stacking sequence is in no way optimised for strength in these regions, as considered for example by Dillinger et al. \cite{dillinger}. In practice, this change of stacking sequence is easy for the user to specify using a \verb=.csv= file specifying different \texttt{Regions} for each segment of the wingbox and providing the required different stacking sequence. The wingbox geometry is again achieved by specifying a \texttt{gridTransformation()}, which now becomes slightly more complex, in order to handle each of the different regions. To create the closed curve of this wingbox, periodic boundary conditions are imposed. In this application, we consider two forms of loading. Firstly, an internal pressure of $0.109$MPa, arising from the fuel, is applied to the internal surface. Secondly, a thermal pre-stress induced by the manufacturing process is imposed, using the user-defined function \texttt{evaluateHeat()}.

\begin{figure}[t]
\centering
\includegraphics[height = 4 cm]{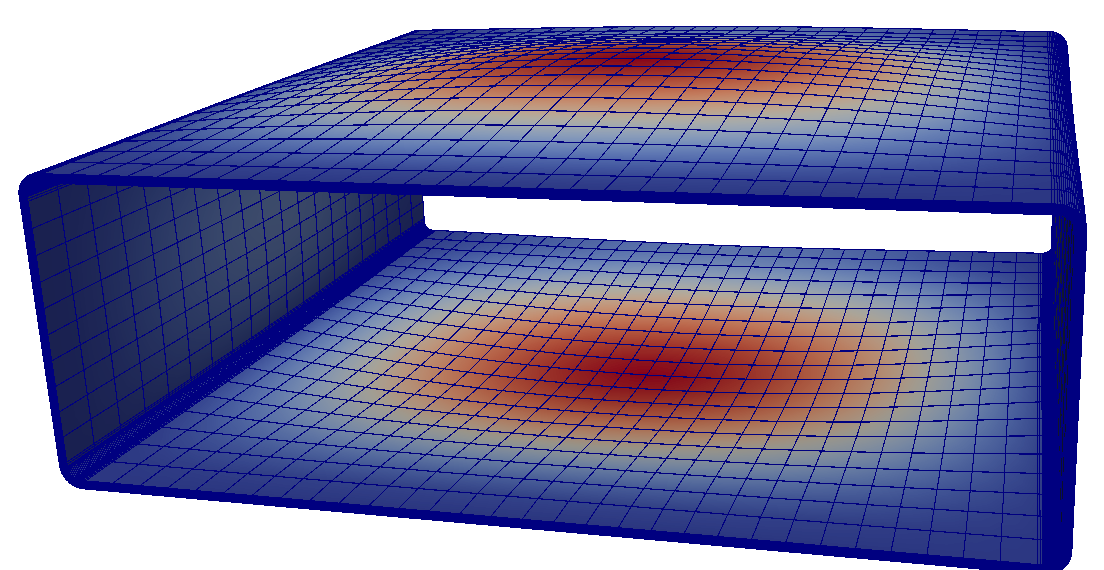} $\quad$
\includegraphics[height = 4 cm]{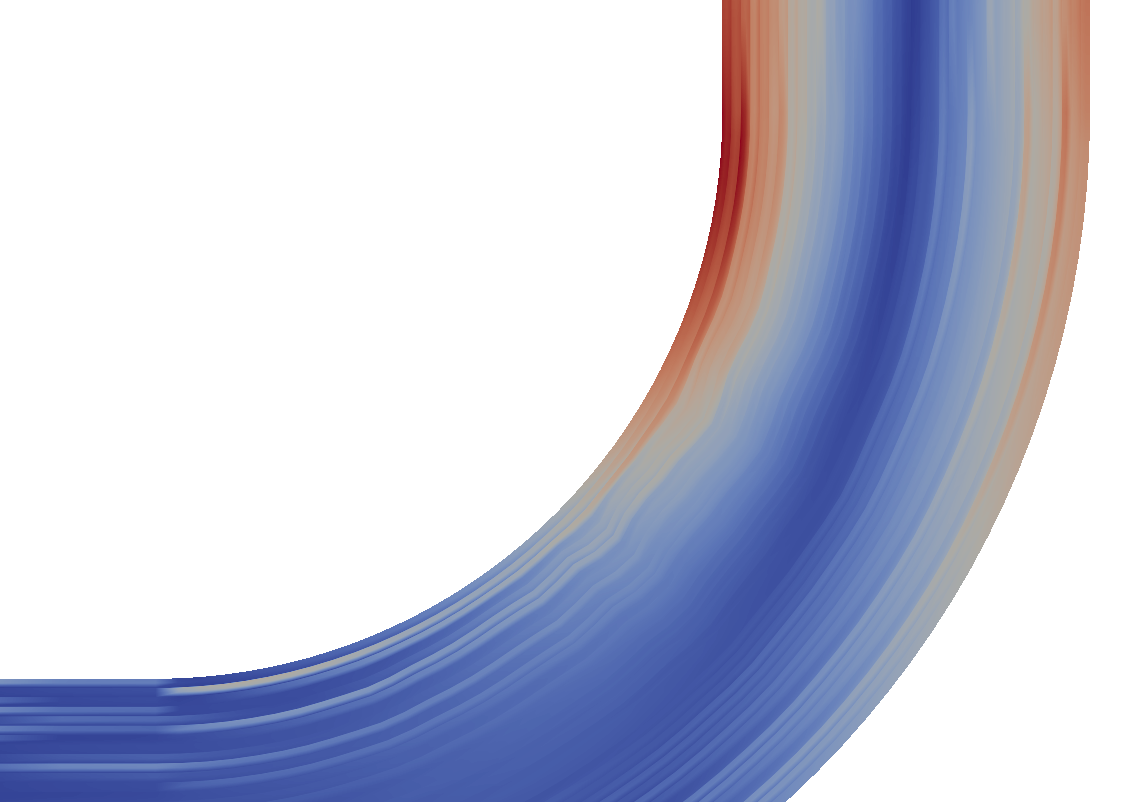}
\caption{\label{fig:wingbox2}FE solution for Example 3: (Left) Overall deformation of the wingbox with colours showing the magnitude of the displacements in cm. (Right) Camanho failure criterion \eqref{eqn:Hashin} in a close-up of the corner containing the wrinkle defect.}
\end{figure}


\smallskip\noindent
We approximate the influence of the ribs that constrain the wingbox in the $y$ direction, by clamping all degrees of freedom at one end, whilst tying all other degrees of freedom at the other end using a multipoint constraint. Elements to be included in the multipoint constraint are marked with the user-defined function \texttt{isMPC(FieldVec\& x)}. A localised wrinkle defect is introduced into one of the corner radii, as shown in Fig.~\ref{fig:wingboxgeom}. The defect is introduced by adapting the function \texttt{gridTransformation()}. The wrinkle geometry is defined by a random field, parameterised by a Karhunen-Lo\'{e}ve expansion. The actual parametrisation of the wrinkle is chosen to match an observed defect in a CT-Scan of a real corner section. Further details of this methodology are provided in Sandhu {\em et al.}~\cite[Sec.~3]{San18}.

\begin{figure}[t]
\centering
\includegraphics[width = 0.42\linewidth]{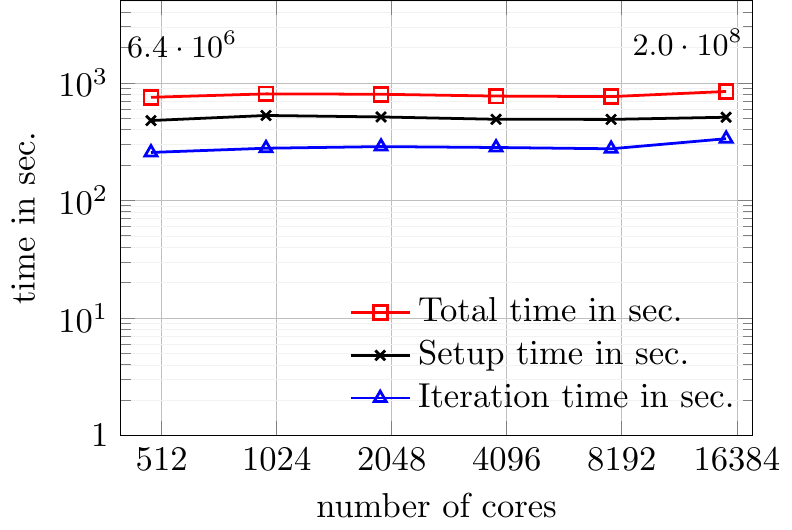}
\includegraphics[width = 0.42\linewidth]{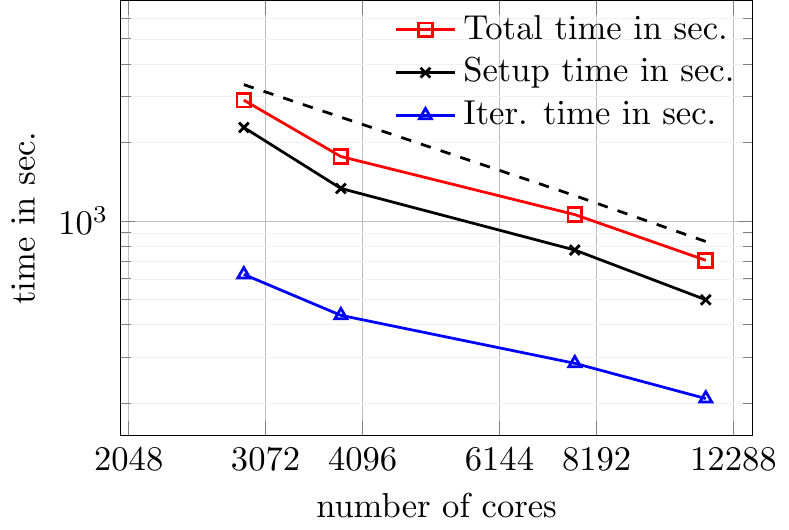}
\caption{\label{fig:performance}Parallel performance of \texttt{dune-composites} on \textsc{Archer}: (Left) A weak scaling test, as summarised in Table \ref{table-wingbox}. (Right) A strong scaling test using Setup~5 in Table~\ref{table-wingbox}, with the dashed line showing perfect scaling, as summarised in Table \ref{table-strongscaling}.}
\end{figure}

\smallskip\noindent
We firstly carry out a weak scaling experiment, increasing the
problem size proportionally to the number of cores used. For iterative solvers that scale optimally with respect to problem size and with respect to the number of cores, the computational time should remain constant. To scale the problem size as the number of cores $N_{cores}$ grows, we refine the mesh, doubling the number of elements as we double the number of cores. The number or elements for each setup are detailed in Table \ref{table-wingbox}, separately listing the number of elements across the spar and the cover, around the corners and along the length of the wingbox. The defective corner, denoted $R_d$, contains twice as many elements as the other three corners, denoted by $R_{nd}$. Table \ref{table-wingbox} also details the resulting number of degrees of freedom, iteration numbers for the preconditioned CG, an estimate of the condition number {\color{black} of the preconditioned system matrix $\bf M_{AS,2}^{-1} A$}, the dimension of the coarse space $\dim V_H$, as well as the total run time. 
Fig.~\ref{fig:performance} (left) shows that the weak scaling of the iterative CG solver in \texttt{dune-composites}  with \textsc{GenEO} preconditioner is indeed almost optimal up to at least 15,360 cores (the limiting capacity available on \textsc{Archer} for our experiments). We also include a more detailed subdivision of the computational time into Setup Time (for the assembly of the FE stiffness matrix and for the construction of the \textsc{GenEO} coarse space) and Iteration Iime (for the preconditioned CG iteration). Both scale almost optimally. This test demonstrates the capability of increasing the size of the tests at a nearly constant run time and thus, to solve a problem with $200$ million degrees of freedom in just over $13$ minutes.

\begin{table}[t]
\centering

\begin{tabular}{|c|c|c|c|c|c|c|c|c|c|c|c|}
\hline
Setup & $N_{cores}$  & Spar & Cover & $R_d$ & $R_{nd}$ &  Length  & DOF & iter. & $\kappa$ & $\dim V_H$ & Time (sec) \\
\hline 
1 & 480  & 34 & 14 & 40 & 20 & 20 &  $6.4\cdot10^6$ & 156 & 445 & 5025 & 734\\ 
2 & 960 & 34 & 14 & 40 & 20 & 40 &  $1.3\cdot10^7$ & 154 & 421 & 7840 & 806\\ 
3 & 1920 & 68 & 28 & 80 & 40 & 40 &  $2.6\cdot10^7$& 152 & 322 & 18752 & 800\\ 
4 & 3840 & 68 & 28 & 80 & 40 & 80 &  $5.1\cdot10^7$& 144 & 287 & 29444 & 772\\ 
5 & 7680 & 216 & 64 & 80 & 40 & 80 &  $1.0\cdot10^8$& 132 & 303 & 50930 & 764\\ 
6 & 15360 & 216 & 64 & 80 & 40 & 160 &  $2.0\cdot10^8$& 102 & 245 & 94527 & 845 \\ 
\hline
\end{tabular}
\caption{\label{table-wingbox}Details of the six setups and results used in the weak scaling test. In all of the tests, we used two layers of 20-node serendipity elements per fibrous layer and only one layer of elements in each of the interface layers. The number of elements per core was fixed at $2808$.}
\end{table}

\smallskip\noindent
Next we carry out a small strong scaling experiment. The mesh is that of Setup~5 in Table~\ref{table-wingbox} and the results of the strong scaling test are given in Figure \ref{fig:performance} (right) and in Table~\ref{table-strongscaling}.  We see that the iterative CG solver in \texttt{dune-composites} with GenEO preconditioner scales almost optimally to at least $11320$ cores, with the time taken approximately halving as the number of cores is doubled. Again, both the Setup and the Iteration scale optimally.


\begin{table}[t]
\begin{tabular}{|c|c|c|c|c|c|c|c|}
\hline
$N_{cores}$  & Elements per Core & $\mbox{dim}(V_H)$ & it. & $T_{it}$ & $T_{setup}$ & $T_{total}$ & Total Core Time (days) \\ \hline
2880     & 3132              & 18843             & 167 & 623      & 2283        & 2906        & 96.9                 \\
3840     & 2340              & 26333             & 153 & 434      & 1332        & 1766        & 78.5           \\
7680     & 2008              & 52622             & 132 & 284      & 773         & 1057        & 94.0               \\
11320    & 1392              & 78233             & 162 & 208      & 498         & 706         & 92.5             \\ \hline
\end{tabular}
\caption{\label{table-strongscaling}Strong scaling test with Setup~5 in Table~\ref{table-wingbox}, demonstrating near optimal strong scaling up to at least $11,320$ cores.}
\end{table}


{\color{black}
\section{Subsurface Flow Application: Strong Scaling for the SPE10 Benchmark}\label{sec:darcy}

\begin{figure}[tb]
 \centering
 \includegraphics[width=0.48\textwidth]{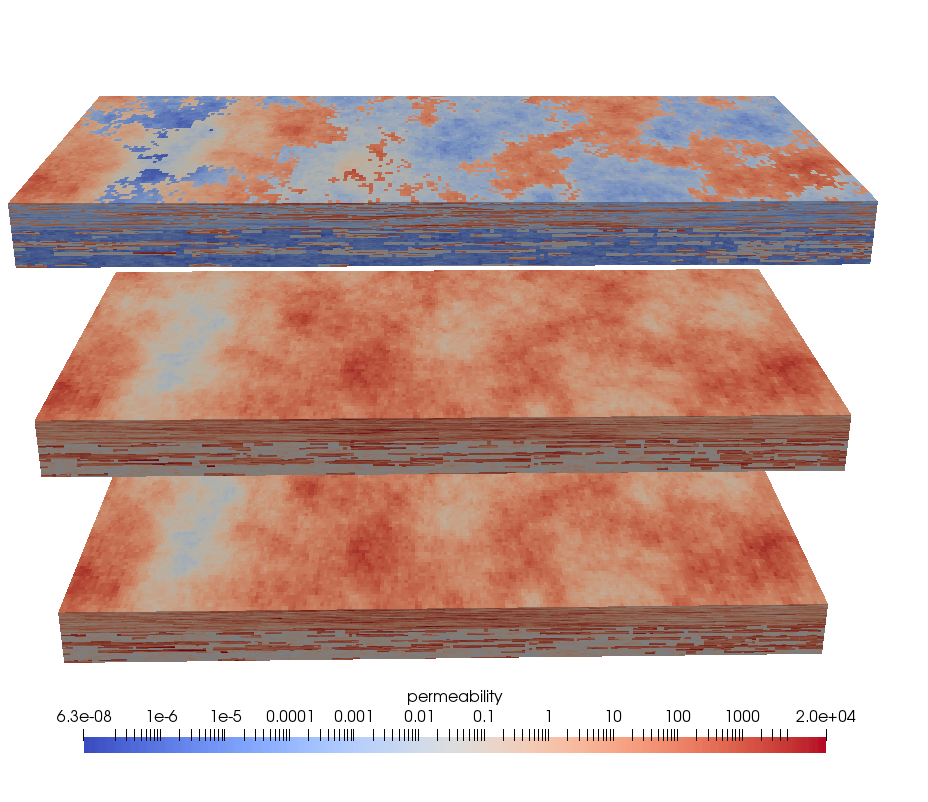}
  \includegraphics[width=0.48\textwidth]{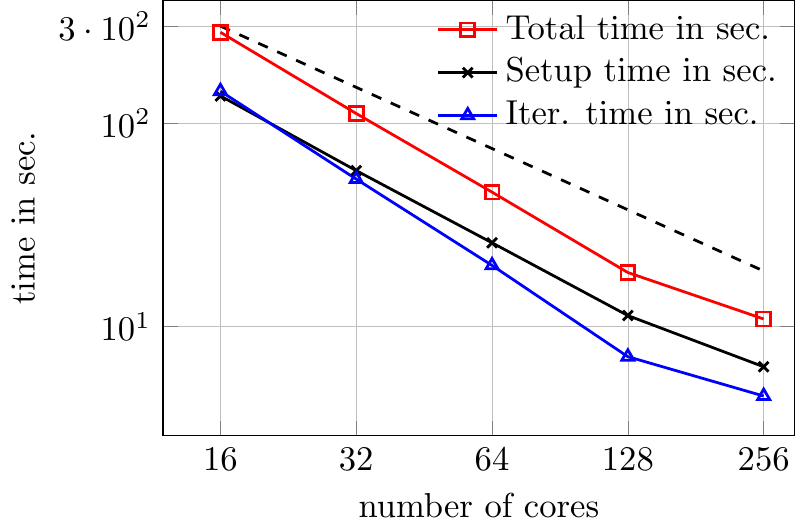}
 \caption{Left: Logarithm of the permeability field $K$ for the SPE10 benchmark, from bottom to top:  $K_x$, $K_y$ and $K_z$. Right: A strong scaling test using the SPE10 dataset, with the dashed line showing perfect scaling. }
 \label{fig-perm}
\end{figure}

\smallskip\noindent In this section, we apply the GenEO solver to an elliptic partial differential problem outside of composites modelling, demonstrating its scalability and robustness on a subsurface flow problem in a highly hetereogeneous medium. A challenging test case in the computational geosciences is the SPE10 benchmark \cite{spe10}. This problem features high contrast, heterogeneous coefficients which challenge most iterative solvers \cite{eike}. 

\smallskip
We consider the SPE10 domain $\Omega := [0, 1200] \times [0, 2200] \times [0, 170]$ (feet), divided into a tensor product grid $\mathcal T_h$ with $60\times220\times85 = 1.122\times 10^6$ cells. The domain $\Omega$ has the boundary $\partial \Omega = \Gamma_D \cup \Gamma_N$, where we define $\Gamma_D := \{{\bf x} \in \partial \Omega : z = 0\}$ as the Dirichlet part of the boundary and ${\bf n} \in \mathbb R^3$ as the outward normal to $\partial\Omega$. We calculate the steady-state fluid pressure $u({\bf x}) \in \Omega$ which obeys Darcy's law. This is given by the linear, scalar elliptic partial differential equation
 \begin{equation}\label{eqn:Darcy}
 -\nabla\cdot({\bf K}({\bf x})\nabla u) = f, \quad \forall {\bf x} \in \Omega
 \end{equation}
subject the boundary conditions
\begin{equation}
u({\bf x}) = 0 \quad \mbox{on} \quad \Gamma_D \quad \mbox{and} \quad  -{\bf K}({\bf x})\nabla u \cdot {\bf n} = 0  \quad \mbox{on} \quad \Gamma_N = \partial\Omega \mbox{\textbackslash} \Gamma_D\,.
\end{equation}
The SPE10 dataset gives a spatially varying permeability tensor
$$
{\bf K}({\bf x}) =
\begin{bmatrix}
K_x({\bf x}) & 0 & 0 \\
0 & K_y({\bf x}) & 0 \\
0 & 0 & K_z({\bf x})
\end{bmatrix} \quad \forall {\bf x} \in \Omega.
$$
Figure \ref{fig-perm}(left) shows the permeability field, it is constant in each cell, but varies strongly over the domain. The parameters $K_x$ and $K_y$ vary from $6.65\times10^{-4}$ to $2.0\times10^4$ and the parameter $K_z$ varies from  $6.65\times10^{-8}$ to $6.0\times 10^3$.
We define the function space for the pressure $u$ as
$V := \{v \in H^1(\Omega) : v({\bf x}) = 0 \;, \ {\bf x} \in \Gamma_D \}$, and choose the finite element space $V_h \subset V$ to be the set of continuous, piecewise linear functions on $\mathcal T_h$. The finite elemen discretisation of \eqref{eqn:Darcy} then reads:  Find $u_h \in V_h$ such that
\begin{equation}\label{eq-fe-darcy}
\int_\Omega {\bf K}({\bf x})\nabla u_h \cdot \nabla v_h \; dx = \int_\Omega f v_h\;dx \quad \forall \mathbf{v}_h \in V_h.
\end{equation}
By defining $u_h = \sum_{i=1}^N u^{(i)}\phi_i({\bf x})$, again we obtain the sparse system of equations 
$$
{\bf A}{\bf u} = {\bf b}, \quad \mbox{where} \quad{\bf u} = [u^{(1)},u^{(2)},\ldots,u^{(N)}]^T$$ 
is a vector pressures at each cell vertex, which is assembled element-wise from \eqref{eq-fe-darcy} using standard Gaussian integration. The  source term in our experiments is assumed to be uniform $f \equiv 1$.

  In Figure \ref{fig-perm}(right) and Table \ref{table-spe10} we show a small strong scaling experiment performed with this challenging setup. The parameter contrast for this benchmark is on the order of $10^{11}$. Nevertheless, we see that the iterative CG solver in \texttt{dune-composites} with GenEO preconditioner scales almost optimally to at least $256$ cores. At $512$ cores with only around $2000$ elements per core the strong scaling begins to break down due to the communication overhead.  Due to the layered structure of the material parameters our domain decomposition is two dimensional. Each subdomain includes the full length in $z$-direction. We used a minimal overlap of only one element. Table \ref{table-spe10} also details the number of cores, size of the coarse space $\dim V_H$, iteration numbers for the preconditioned CG, as well as the time spent in CG iterations, setup time and total run time.
  
    \begin{table}[tb]
 \begin{tabular}{|c|c|c|c|c|c|}
 \hline
$N_{cores}$  & $\mbox{dim}(V_H)$ & it. & $T_{it}$ & $T_{setup}$ & $T_{total}$ \\
       \hline
16  &  149  &   167 & 136.11  & 143.621  & 279.721 \\
32  &  225  &	203 & 58.42	 & 53.065	& 111.485\\
64	&  379  &	206	& 25.81 & 19.982	& 45.787\\
128	&  527  &	224	& 11.32	 & 7.107	& 18.427\\
256	&  930  &	232	& 6.34	 & 4.552	& 10.892\\
512	&  1737 &	234	& 5.18 & 3.795	& 8.975\\
 \hline
 \end{tabular}
 \caption{A strong scaling test using the SPE10 dataset.}
 \label{table-spe10}
 \end{table}
  
}

\section{Discussion \& Future Developments} \label{sec:conclusions}

\noindent
In this paper, we describe the new high performance package \texttt{dune-composites}, designed to solve massive finite element problems for the anisotropic linear elasticity equations. The paper provides both the mathematical foundations of the methods, their implementation within a state-of-the-art software platform on modern distributed memory computer architectures, as well as details of how to set up a problem and carry out an analysis, illustrated via a series of increasingly complex examples. In addition, we demonstrate the scalability of the new solver on over $15,000$ cores on the UK national supercomputer \textsc{Archer}, solving industrially motivated problems with over $200$ million degrees of freedom within minutes. This scale of computations brings composites problems that would otherwise be unthinkable into the feasible range.

\smallskip
\noindent


\smallskip\noindent
The disadvantage of \texttt{dune-composites} as a package over commercial counterparts is currently the limited functionality in considering more general problems; this includes complex geometries, unstructured grids and nonlinear problems. This is the first release of \texttt{dune-composites}, and therefore the functionality is naturally still limited, but it will increase over time, driven by the industrial questions we seek to solve as a community of developers. There are currently four key areas of active development:
\begin{itemize}
\item {\bf Multiscale Methods}: There is a strong connection between coarse spaces for domain decomposition methods, as developed and implemented within the \textsc{GenEO} preconditioner, and multiscale discretisation methods, such as generalised multiscale FEs (GMsFE) \cite{Efe13}. In fact,  the \textsc{GenEO} coarse space provides and natural multiscale method. Current research \cite{Dod17} is aiming to provide this functionality as an embedded solution scheme within \texttt{dune-composites}.

\item {\bf Nonlinear Mechanics}: Currently the analysis implemented in the package is linear (static and thermal). Naturally, modelling failure propagation in composites is a nonlinear problem. Current research in implementing nonlinear material models includes cohesive zone models and ply damage models, for integration within the existing framework.

{\color{black}

\item {\bf Uncertainty Quantification}: A significant motivation for developing efficient robust solvers is to enable stochastic studies in which many simulations of a model with different parameters are required. The package has already been used to study the effects of wrinkle defects in composite strength \cite{San18}, and ongoing research is exploring Bayesian parameter estimation using multilevel Markov chain Monte Carlo methods~\cite{Dod15SIAM, Dod19SIAM}.}

\item {\bf GUI Development}: A current limitation of the package is that its application to new models or geometries requires a basic knowledge of \texttt{C++}  and command line programming. Active development with software engineers is seeking to provide a simple Graphical User Interface (GUI) to enable application-focused research with \texttt{dune-composites} without extensive programming experience.

\end{itemize}

 \section{Acknowledgements}
 This work was supported by an EPSRC Maths for Manufacturing grant (EP/K031368/1).
 Richard Butler holds a Royal Academy of Engineering-GKN Aerospace Research Chair in
 Composites. This work used the ARCHER UK National Supercomputing Service (\url{http://www.archer.ac.uk}).

\section*{References}
\bibliographystyle{model1-num-names}
\bibliography{references.bib}

\end{document}